\documentclass[11pt]{article}
\usepackage{graphicx,amsthm,graphicx,fancyhdr,mathrsfs}
\usepackage{amsfonts}
\usepackage{amsmath}
\usepackage{amssymb}
\usepackage{url}
 \usepackage[colorlinks=true,citecolor=blue]{hyperref}
\usepackage{fancyhdr}
\usepackage{indentfirst}
\usepackage{enumerate}
\usepackage{amsthm}
\usepackage{color}
\usepackage{comment}
\usepackage{dsfont}
\usepackage{natbib}
\usepackage[misc]{ifsym}
\usepackage[toc,page]{appendix}
\usepackage{multirow}
\usepackage{float}
\usepackage{bm}
\usepackage{booktabs}
\usepackage{todonotes}

\def\d{\,\mathrm{d}}
\def\laweq{\buildrel \mathrm{d} \over =}

\newcommand{\VaR}{\mathrm{VaR}}

\newcommand{\ES}{\mathrm{ES}}

\renewcommand{\H}{\mathcal{H}}

\newcommand{\I}{\mathcal{I}}

\newcommand{\X}{\mathcal{X}}

\newcommand{\LE}{\mathrm{ES}^{L}}

\newcommand{\E}{\mathbb{E}}
\newcommand{\R}{\mathbb{R}}
\newcommand{\N}{\mathbb{N}}
\renewcommand{\L}{\mathcal{L}}
\newcommand{\M}{\mathcal{M}}
\newcommand{\p}{\mathbb{P}}
\renewcommand{\P}{\mathbb{P}}

\newcommand{\id}{\mathds{1}}

\renewcommand{\ge}{\geqslant}
\renewcommand{\le}{\leqslant}

\renewcommand{\leq}{\leqslant}
\renewcommand{\epsilon}{\varepsilon}

\theoremstyle{plain}
\newtheorem{theorem}{Theorem}
\newtheorem{corollary}{Corollary}
\newtheorem{lemma}{Lemma}
\newtheorem{proposition}{Proposition}
\theoremstyle{definition}
\newtheorem{definition}{Definition}
\newtheorem{example}{Example}

\newtheorem{remark}{Remark}

\DeclareMathOperator*{\argmax}{arg\,max}
\DeclareMathOperator*{\argmin}{arg\,min}

\usepackage{tikz}

\renewcommand{\cite}{\citet}

\makeatletter
\DeclareRobustCommand{\bsquare}{%
	\mathop{\vphantom{\sum}\mathpalette\bigstar@\relax}\slimits@
}

\newcommand{\bigstar@}[2]{%
	\vcenter{%
		\sbox\z@{$#1\sum$}%
		\hbox{\resizebox{.9\dimexpr\ht\z@+\dp\z@}{!}{$\m@th\dsquare$}}%
	}%
}
\makeatother

\usepackage[onehalfspacing]{setspace}

\newcommand{\dsquare}{\mathop{  \square} \displaylimits}

\begin{document}

	\title{Optimizing distortion riskmetrics with distributional uncertainty}

	\author{Silvana M. Pesenti\thanks{Department of Statistical Sciences, University of Toronto, Canada. \Letter~\texttt{silvana.pesenti@utoronto.ca}} \and Qiuqi Wang\thanks{Department of Statistics and Actuarial Science, University of Waterloo, Canada. \Letter~\texttt{q428wang@uwaterloo.ca}}\and Ruodu Wang\thanks{Department of Statistics and Actuarial Science, University of Waterloo, Canada. \Letter~\texttt{wang@uwaterloo.ca}}}
	\maketitle 
	
	\begin{abstract}
		Optimization of distortion riskmetrics with distributional uncertainty has wide applications in finance and operations research. Distortion riskmetrics include many commonly applied risk measures and deviation measures, which are not necessarily monotone or convex. 
		One of our central findings is  a unifying result that allows us to convert an optimization of a non-convex distortion riskmetric with distributional uncertainty to a convex one, leading to great tractability. 
		A sufficient condition to the unifying equivalence result is the novel notion of closedness under concentration, a variation of which is also shown to be necessary   for the equivalence.
		Our results include many special cases that are well studied in the optimization literature, including but not limited to optimizing probabilities, Value-at-Risk, Expected Shortfall,  Yaari's dual utility, and differences between distortion risk measures, under various forms of distributional uncertainty.
		We illustrate our theoretical results via applications to portfolio optimization, optimization under moment constraints, and preference robust optimization.
		
		~
		
	\end{abstract} 
	
	~
	
	\section{Introduction}

	Riskmetrics, such as measures of risk and variability, are common tools to represent preferences, model decisions under risks, and quantify different types of risks.  
	To fix terms, we refer to \emph{riskmetrics} as any mapping from a set of random variables to the real line, and \emph{risk measures} as riskmetrics that are monotone in the sense of \cite{ADEH99}.

	In this paper, we focus on \emph{distortion riskmetrics} which is a large class of commonly used measures of risk and variability; see \cite{WWW20} for the terminology ``distortion riskmetrics".
	Distortion riskmetrics include
	L-functionals \citep{HR09} in statistics, Yaari's dual utilities \citep{Y87} in decision theory, distorted premium principles \citep{WYP97} in insurance, and  spectral  risk measures \citep{A02} in finance; see \cite{WWW20} for further examples. 
	After a normalization, increasing distortion riskmetrics are  \emph{distortion risk measures}, which include,
	in particular, the two most important risk measures used in current banking and insurance regulation, the Value-at-Risk (VaR) and the Expected Shortfall (ES). 
	Moreover, convex distortion riskmetrics are the building blocks (via taking a supremum) for all convex risk functionals \citep{LCLW19},
	including   classic risk measures \citep{ADEH99,FS02} and deviation measures \citep{RUZ06}. 
	
    When riskmetrics are evaluated on distributions that are subject to uncertainty, decisions should be taken with respect to the worst (or best) possible values a   riskmetric  attains over a set of alternative distributions; giving rise to the active subfield of distributionally robust optimization. 
    The set of alternative distributions, the \emph{uncertainty set}, may be characterized by moment constraints (e.g., \cite{P07}),   parameter uncertainty (e.g., \cite{DY10}),   probability constraints (e.g., \cite{WKS14}), and distributional distances (e.g., \cite{BM19}),  amongst others. 
    Popular distortion risk measures such as VaR and ES are studied extensively in this context; see e.g., \cite{NPS08} and \cite{ZF09}.

    Optimization of  convex distortion risk measures, i.e., distortion riskmetrics with an increasing and concave distortion function,  is relatively well understood under distributional uncertainty; see  
        \cite{CRV18}, \cite{L18},  and \cite{LCLW19} for some recent work.
		 Nevertheless, many distortion riskmetrics are not convex or monotone. For example, in the Cumulative Prospect Theory of \cite{TK92}, the distortion function is typically assumed to be inverse-S-shaped; in financial risk management, the popular risk measure $\VaR$ has a non-concave distortion function, and the inter-quantile difference \citep{WWW19} has a distortion function that is neither concave nor monotone. Another example is the difference between two distortion risk measures, which is clearly not increasing or convex in general. 
	Optimizing  non-convex distortion riskmetrics under distributional uncertainty is difficult and results are available only for special cases; see \cite{LSWY18},  \cite{CLW18}, \cite{ZS18}, \cite{WXZ19}, and  \cite{BPV20}, all with an increasing distortion function.

    There is, however, a notable common feature in  the above mentioned literature when a non-convex distortion risk metric is involved. For numerous special cases, one often obtains an equivalence between the optimization problem with non-convex distortion riskmetric and that with a convex one.
	Inspired by this observation, the aim of this paper is to address: 
	\begin{center} What conditions provide equivalence between a non-convex riskmetric and a convex one in the setting of distributional uncertainty?
	\end{center}
	An answer to this question is still missing in the literature.  In this sense, we offer a novel perspective on distributionally robust optimization problems  by converting non-convex optimization problems to their convex counterpart. 
	    Transforming a non-convex to a convex optimization problem through approximation and via a direct equivalence has   been studied by \cite{ZKR13} and \cite{CLM20}. Both contributions, however, consider uncertainty sets described by some special forms of constraints. A unifying framework applicable to numerous uncertainty sets and the entire class of distortion riskmetrics is however missing and at the core of this paper.
	

	The main novelty of our results is three-fold: first, we obtain a unifying result (Theorem \ref{th:1}) that allows, under distributional uncertainty, to convert an optimization problem of a non-convex distortion riskmetric  to an optimization problem with a convex one. The result covers, to the authors' best knowledge, all known equivalences between optimization problems of non-convex and convex riskmetrics with distributional uncertainty. 
    The proof requires  techniques beyond the ones used in the existing literature, as we do not make assumptions such as monotonicity, positiveness, and continuity.  
    Our framework can also be easily applied to settings with atomic probability space or with uncertainty sets of multi-dimensional distributions.
	Second, we introduce the concept of \emph{closedness under concentration} as a sufficient condition to establish the equivalence, and it is also a necessary condition on the set of optimizers given that the equivalence holds (Theorem \ref{th:nec}). We show how the properties of closedness under concentration within a collection of intervals $\I$ and closedness under concentration for all intervals can easily be verified and provide numerous examples. Third, the classes of distortion riskmetrics  and uncertainty formulations  considered in this paper include all special cases studied in the literature; examples are presented in Sections \ref{sec:main}-\ref{sec:multi}. In particular,  
our class of riskmetrics include  all practically used risk measures and variability measures (some via taking a sup),
dual utilities with inverse-S-shaped distortion functions of \cite{TK92}, and 
differences between two dual utilities or distortion risk measures.
 Our uncertainty formulations include both supremum and infimum problems,\footnote{Thus we provide a universal treatment of worst-case and best-case risk values. Calculating best-case risk values allows us to solve economic decision making problems where optimal distributions are chosen to minimize the risk.} moment constraints, convex order/risk measure constraints,  marginal constraints in risk aggregation with dependence uncertainty (e.g., \cite{EWW15}), preference robust optimization (e.g., \cite{AD15}  and \cite{GX20}), and some one-dimensional and multi-dimensional uncertainty sets induced by Wasserstein metrics.

	The great generality distinguishes our work from the large literature on distributional robust optimization cited above.  
	Our work is of analytical and probabilistic nature, and we focus on   theoretical equivalence results which will be also illustrated via numerical implementations.
	The target problems are formulated in Section \ref{sec:2}.
	Sections \ref{sec:main} is devoted to our main contribution of the equivalence of non-convex and convex optimization problems with distributional uncertainty.
	We illustrate by many examples the concepts of closedness under conditional expectation and closedness under concentration, and distinguish them in several practical settings.
	Section \ref{sec:multi} demonstrates the equivalence results in multi-dimensional settings. In addition to a general multi-dimensional model with a concave loss function,  we solve a robust risk aggregation problem with ambiguity on both the marginal distributions and the dependence structure.
	In Section \ref{sec:two_con},  our results  are used to solve optimization problems with uncertainty sets defined via moment constraints. In particular, we generalize a few well-known results in the literature on optimization and worst-case values of risk measures. 
	Sections \ref{sec:app} and \ref{sec:app_num} contain numerical illustrations of optimizing differences between two distortion riskmetrics, portfolio optimization, and  preference robust optimization.  Some concluding remarks are put in Section \ref{sec:conclude}. Proofs of all results are relegated to Appendix \ref{app:proofs}.

	\section{Distortion riskmetrics with distributional uncertainty}
	\label{sec:2}
	
	\subsection{Problem formulation}
	Throughout, we work with an atomless probability space $(\Omega,\mathscr F,\p)$. For $n\in\N$,   $A$ represents a set of actions,  $\rho$ is an objective functional,  
	$f: A\times\R^{n} \to \R$ is a loss function, 
	and  $\mathbf{X}$ is an $n$-dimensional random vector with distributional uncertainty.
	Many problems in distributionally robust optimization have the form 
	\begin{equation}\label{eq:opt}
		\min_{\mathbf{a}\in A}\, \sup_{F_\mathbf{X}\in \widetilde{\mathcal M} }\; \rho(f(\mathbf{a},\mathbf{X})),
	\end{equation}
	where $F_\mathbf{X}$ denotes the distribution of $\mathbf{X}$ and $\widetilde{\mathcal M}$ is a set of plausible distributions for $\mathbf{X}$. 
	We will first focus on the inner problem
	\begin{equation}\label{eq:inner}
		\sup_{F_\mathbf{X}\in \widetilde{\mathcal M} }\; \rho(f(\mathbf{a},\mathbf{X})),
	\end{equation}
	which we may rewrite as 
	\begin{equation}\label{eq:yproblem}
		\sup_{F_Y\in \mathcal M } \rho(Y),
	\end{equation}
	where $F_Y$ denotes the distribution of $Y$ and $\mathcal M$ is a set of distributions on $\R$. We suppress the reliance on $\mathbf{a}$ as it remains constant in the inner problem \eqref{eq:inner}.  The supremum in \eqref{eq:yproblem} is typically referred to as the \emph{worst-case risk measure} in the literature if $\rho$ is monotone.\footnote{A risk measure $\rho:\L^p\to\R$ is \emph{monotone} if $\rho(X)\le \rho(Y)$ for all $X,Y\in\L^p$ with $X\le Y$.}
The problem \eqref{eq:yproblem} can also represent an optimal decision problem, where $\rho$ is an objective to maximize, and  a decision maker chooses an optimal distribution from the set $\M$ which is interpreted as an action set instead of an uncertainty set (i.e., no uncertainty in this problem).
	Since the two problems share the same mathematical formulation \eqref{eq:yproblem}, we will  navigate through our results mainly with the first interpretation of worst-case risk under uncertainty. 
	
	We denote by $\L^p$, $p\in [1,\infty)$, the space of random variables with finite $p$-th moment. Let $\L^\infty$ represent the set of bounded random variables and let $\L^0$ represent the space of all random variables.
	Denote by $\mathcal H$ the set of  functions {$h:[0,1]\mapsto\R$} with bounded variation satisfying 
	$h(0)=0$. 
	For  $p\in [1,\infty]$ and $h\in\H$, a \emph{distortion riskmetric} $\rho_h: \mathcal  L^p\to \R$ is defined as 
	\begin{equation}\label{eq:disrep}
		\rho_{h}(Y) = \int_0^\infty  h (\p(Y>x))  \d x + \int_{-\infty}^0 (h (\p(Y> x)) -h(1))\d x,~~Y\in \mathcal L^p,
	\end{equation}
	whenever the above integrals are finite; see Proposition \ref{co:finiteness} below for a sufficient condition.   The function $h\in\mathcal H$ is called a \emph{distortion function}.
	Note that we allow $h$ to be non-monotone; if $h$ is increasing and $h(1)=1$, then $\rho_h$ is a distortion risk measure. The distortion riskmetric $\rho_h$ is convex if and only if $h$ is concave; see \cite{WWW19} for this and other properties of $\rho_h$.
	
	In this paper  we consider the objective functional $\rho$ in \eqref{eq:opt} to be a distortion riskmetric $\rho_h$ for some $h\in \mathcal H$, as the class of distortion riskmetrics  includes a large class of objective functionals of interest. 
	Note that a general analysis of \eqref{eq:yproblem}  also covers the infimum problem
	$ \inf_{F_Y\in \mathcal M} \rho_h(Y)$, 
	since $-\rho_h=\rho_{-h}$ is again a distortion riskmetric. This illustrates an advantage of studying distortion riskmetrics over monotone ones, as our analysis unifies best- and worst-case risk evaluations.  Best-case risk measures are also of practical importance. 
	In particular, they may represent risk minimization problems through the second interpretation of \eqref{eq:yproblem}, where $\M$ represents a set of possible actions (see Section \ref{sec:example} for some examples).

	
	If $\rho_h$ is not convex, or equivalently, $h$ is not concave,
	optimization problems of the type \eqref{eq:yproblem} are often highly nontrivial.  
	However, 
	the optimization problem of  maximizing $
		  \rho_{h^*}(Y)$ over ${F_Y\in \mathcal M }$,
	where $h^*$ is the smallest concave distortion function dominating $h$, is convex and can often be solved relatively easily either analytically or through numerical methods. 
	Clearly, since $\rho_{h^*}\ge \rho_h$, we have
	$$
	\sup_{F_Y\in \mathcal M } \rho_h(Y) \le \sup_{F_Y\in \mathcal M }  \rho_{h^*}(Y)
	$$
	and one naturally wonders when the above inequality holds as equality, that is, under what conditions, it holds 
	\begin{equation}\label{eq:main}
		\sup_{F_Y\in \mathcal M } \rho_h(Y) = \sup_{F_Y\in \mathcal M }  \rho_{h^*}(Y)\,.
	\end{equation}
	The main contribution of this paper is a sufficient condition on the uncertainty set $\mathcal M$ that guarantees equivalence of these optimization problems, that is \eqref{eq:main} holds.  We will also obtain a necessary condition.
	If \eqref{eq:main} holds, then the non-convex problem (the left-hand side of    \eqref{eq:main}) is converted into the convex problem (the right-hand side of   \eqref{eq:main}), 
	providing huge convenience, which in turn helps to solve the minimax problem \eqref{eq:opt}.


	\subsection{Notation and preliminaries}\label{sec:22}
\begin{figure}[htbp]
		\begin{center}
			\includegraphics[width=0.45\textwidth]{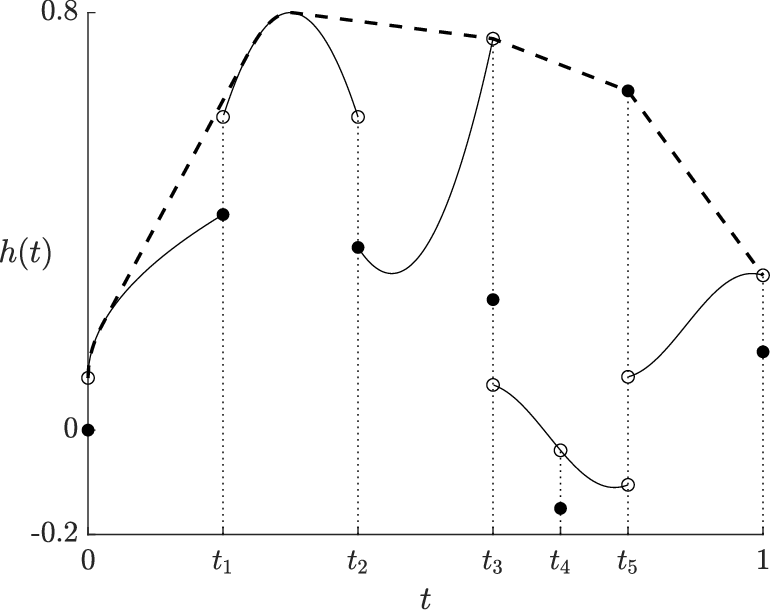}
			 ~~~~ %
			\includegraphics[width=0.45\textwidth]{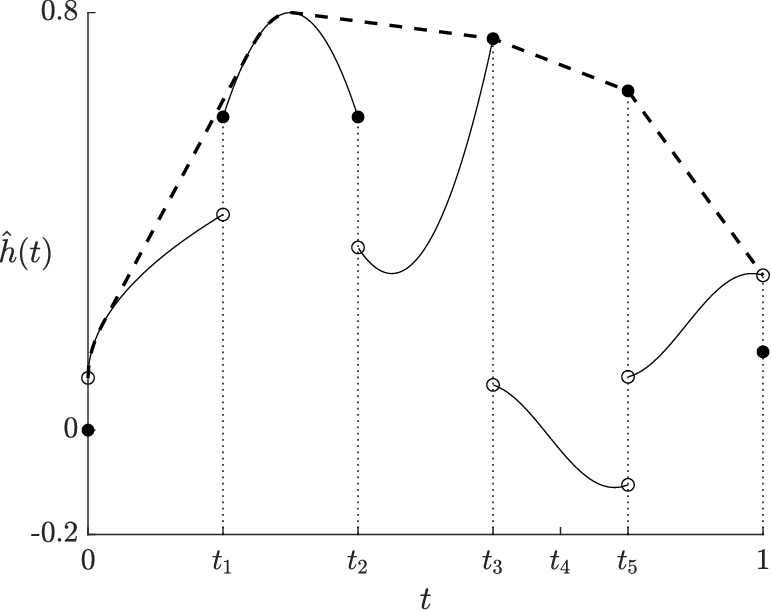}
		\end{center}
		\caption{An example of  $h$ (left) and $\hat{h}$ (right) with the set of discontinuity points $\{t_1,t_2,t_3,t_4,t_5\}$ excluding $0$ and $1$; the dashed lines represent $h^*$ and $(\hat{h})^*$, which are identical by Proposition \ref{lem:extension}}
		\label{fig:hat-h}
		\end{figure}

	For $p\ge 1$ and
	$n\in\N$, we denote by $\M^n_p$ the set of all distributions on $\R^n$ with finite $p$-th moment. Let $\M^n_\infty$ be the set of $n$-dimensional distributions of bounded random variables. For $p\in[1,\infty]$, write $\M^1_p=\M_p$ for simplicity.
	The set inclusion $\subset$ and terms like ``increasing" and ``decreasing" are in the non-strict sense. 
	For $X,Y\in\L^p$, we write $X\laweq Y$ to represent that $X$ and $Y$ have the same distribution.
	For a distribution $F\in \mathcal M_1$, let its left- and right-quantile functions be given respectively by 
	$$F^{-1}(\alpha) = \inf\:\{x\in \R: F(x) \ge \alpha\}~~\text{and}~~F^{-1+}(\alpha) = \inf\:\{x\in \R: F(x) > \alpha\},~~\alpha\in [0,1],$$
	with the convention $\inf(\varnothing)=\infty$. For $x,y\in\R$, we write $x\vee y=\max\{x,y\}$ and $x\wedge y=\min\{x,y\}$. Since $h\in\H$ is of bounded variation, its discontinuity points are at most countable and the left- and right-limits exist at each of these points. We write
	$$h(t^+)=\left\{\begin{array}{l l}
	\lim_{x\downarrow t}h(x), & t\in[0,1),\\
	h(1), & t=1,
	\end{array}\right.~~\text{and}~~ h(t^-)=\left\{\begin{array}{l l}
	\lim_{x\uparrow t}h(x), & t\in(0,1],\\
	h(0), & t=0,
	\end{array}\right.$$ 
	and the upper semicontinuous modification of $h$  is denoted by
 $$ 
			\hat{h}(t)= 
				h(t^+)\vee h(t^-)\vee h(t),~~t\in (0,1),~~\mbox{with }\hat h(0)=0~\mbox{and}~\hat h(1)=h(1).
 $$
		Note that $\hat h(t)=h(t)$ at all continuous points of $h$, and we do not make any modification at the points $0$ and $1$  even if $h$ has a jump at these points. 
	For $h\in \mathcal H$ and $t\in [0,1]$, define its concave and convex envelopes $h^*$ and $h_*$ respectively by
	$$
	h^*(t)=\inf\left\{g(t) : ~g\in \mathcal H,~g\ge h,~g \textrm{ is   concave on}~[0,1]\right\},
	$$
	$$
	h_*(t)=\sup\left\{g(t) : ~g\in \mathcal H,~g\leq h,~g \textrm{ is   convex on}~[0,1]\right\}.
	$$
	Both $h^*$ and $h_*$ are continuous functions  on $(0,1) $ for all $h\in \mathcal H$,
	and if $h$ is continuous at $0$ and $1$, then so are $h^*$ and $h_*$ (see Figure \ref{fig:distortion} below for an illustration of $h$ and $h^*$).
	Denote by $\mathcal H^*$ (resp.~$\mathcal H_*$) the set of   concave (resp.~convex) functions in $\mathcal H$.
	Note that for all $h\in \mathcal H$, we have $h^*\in \mathcal H^*$ and $h_*\in \mathcal H_*$. As a well-known property of the convex and concave envelopes of a   continuous $h$ (e.g., \cite{BC94}),
	 $h^*$ (resp.~$h_*$) differs from $h$ on a union of disjoint open intervals, and $h^*$ (resp.~$h_*$) is linear on these intervals. 
The functions $h$,  $\hat h$, $h^*$ and $(\hat h)^*$  are illustrated in Figure \ref{fig:hat-h}. 

	While in general $\rho_h$ and $\rho_{\hat h}$ are different functionals, one has $\rho_h(Y) = \rho_{\hat h}(Y)$ for any random variable $Y$ with continuous quantile function; see Lemma 1 of \cite{WWW20}.
	Moreover, $h^*=(\hat h)^*\ge \hat h \ge h$ and the four functions are all equal if $h$ is concave.  
	Below, we provide a new result on convex envelopes of distortion functions $h$ that are not necessarily monotone or continuous,  which may be of independent interest.

	  \begin{proposition}\label{lem:extension}
		For any $h\in\H$,  we have $h^*=(\hat h)^*$ and
		the set $\{t\in [0,1]: \hat{h}(t)\ne h^*(t)\}$ is the union of some disjoint open intervals. Moreover,
		$h^*$ is linear on each of the above intervals.
	\end{proposition}
	
	 In the sequel, we mainly focus on $h^*$, which will be useful when optimizing $\rho_h$ in \eqref{eq:yproblem}. A similar result to Proposition \ref{lem:extension} holds for $h_*$, useful in the corresponding infimum problem, where the upper semicontinuous modification of $h$ is replaced by the lower semicontinuous one. This follows directly from Proposition \ref{lem:extension} by setting $g=-h$ which gives $\rho_g=-\rho_h$ and $h_*=-g^*$.
	
	For all distortion functions $h\in \mathcal H$, from Proposition \ref{lem:extension}, there exist (countably many) disjoint open intervals on which $\hat{h}\ne h^*$.
	Using a similar notation to \cite{WXZ19}, we define the set
	$$\mathcal{I}_h=\{(1-b,1-a):\hat{h}\neq h^* \text{ on }(a,b),~\hat{h}(a)=h^*(a),~\hat{h}(b)=h^*(b)\}\,.$$
	The set $\I_h$ is easy to identify in practice; see Section \ref{sec:riskmetrics} for examples of commonly used distortion riskmetrics and their corresponding sets $\I_h$.
	
	\section{Equivalence between non-convex and convex riskmetrics}
	\label{sec:main}

	\subsection{Concentration and the main equivalence result}
	\label{sec:31}
\begin{figure}[t]     
         \centering 
\begin{tikzpicture}
\draw[<->] (0,3.88) -- (0,0) -- (3.4,0); 
\draw[gray,dotted] (3,0) -- (3,3.65); 
\node[below] at (3,0) {$1$};  
  \draw[thick,black,domain=0:3,smooth] plot (\x,{0.4*\x^2});
\node[below] at (0,0) {$0$};    
\node[below] at (1,0) {$\frac 13$};  
\draw[gray,dotted] (1,0) -- (1,0.4); 
\node[below] at (1.5,0) {$\frac 12$};  
\draw[gray,dotted] (1.5,0) -- (1.5,0.9); 
\node[below] at (2,0) {$\frac 23$};  
\draw[gray,dotted] (2,0) -- (2,1.6); 
\draw[red, very thick] (0,0) -- (1,0); 
\draw[red, very thick] (1.5,0) -- (2,0);  
\node[right] at (0,3.9) {$F^{-1}$};  
\end{tikzpicture} ~~~~~~~~~
\begin{tikzpicture}
\draw[<->] (0,3.88) -- (0,0) -- (3.4,0); 
\draw[gray,dotted] (3,0) -- (3,3.65); 
\node[below] at (3,0) {$1$};    
  \draw[thick,black,domain=1:1.5,smooth] plot (\x,{0.4*\x^2});   
\draw[black, thick] (0,0.16) -- (1,0.16); 
  \draw[thick,black,domain=2:3,smooth] plot (\x,{0.4*\x^2});
\draw[black, thick] (1.5,1.2) -- (2,1.2);  
\node[below] at (0,0) {$0$};    
\node[below] at (1,0) {$\frac 13$};  
\draw[gray,dotted] (1,0) -- (1,0.4); 
\node[below] at (1.5,0) {$\frac 12$};  
\draw[gray,dotted] (1.5,0) -- (1.5,0.9); 
\node[below] at (2,0) {$\frac 23$};  
\draw[gray,dotted] (2,0) -- (2,1.6); 
\draw[red, very thick] (0,0) -- (1,0); 
\draw[red, very thick] (1.5,0) -- (2,0);  
\node[right] at (0,3.9) {$(F^{\mathcal I})^{-1}$};  
\end{tikzpicture}
 
\caption{Left panel:  quantile function of $F$; right panel: quantile function of $F^{\mathcal I}$ where $\mathcal I=\{(0,1/3),(1/2,2/3)\}$}
\label{fig:r1-new}
\end{figure}
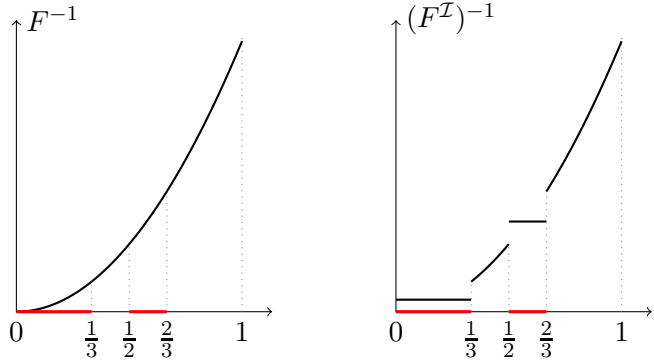

	In this section, we introduce the concept of concentration, and use this concept to explain our main equivalence results, Theorems \ref{th:1} and \ref{th:nec}.  
	For a distribution $F \in \mathcal M_1$
	and an interval $C\subset [0,1]$ (when speaking of an interval in $[0,1]$, we exclude singletons or empty sets),  
	we   define the \emph{$C$-concentration of $F$},
	denote by $F^C$, as the distribution of the random variable
	\begin{equation}
		\label{eq:con}
		F^{-1}(U)\id_{\{U\not \in C\}} +  \E[F^{-1}(U)|U\in C]\id_{\{U  \in C\}},
	\end{equation}
	where $U\sim \mathrm{U}[0,1]$ is a standard uniform random variable. 
	In other words,  $F^C$ is obtained by concentrating the probability mass of $F^{-1}(U)$ on $\{U\in C\}$ at its conditional expectation, whereas the rest of the distribution remains unchanged. 
		For $F \in \mathcal M_1$ and $0\le a < b \le 1$, it is clear that
		the left-quantile function   of 
		$F^{(a,b)}$ is given by
		\begin{equation}\label{eq:quantile}
		F^{-1}(t)\id_{\{t\not \in (a,b]\}}  + \frac{\int_a^b F^{-1}(u) \d u}{b-a} \id_{\{t  \in (a,b]\}},~~t\in[0,1].
		\end{equation}
	For a collection $\I$ of (possibly infinitely many) non-overlapping intervals in $[0,1]$, 
	let $F^{\I}$ be the distribution corresponding to the left-quantile function given by the left-continuous version of
	\begin{equation}\label{eq:invG3}
		F^{-1}(t)\id_{\{t\not \in  \bigcup_{C\in \mathcal I} C\}}  +  \sum_{ C\in \mathcal I}\frac{\int_C F^{-1}(u) \d u}{\lambda(C)} \id_{\{t  \in C \}},~~~~t\in [0,1],
	\end{equation} 
	where $\lambda$ is the Lebesgue measure; see Figure \ref{fig:r1-new} for an illustration.	

\begin{definition}\label{def:con}
Let $\mathcal M$ be a set of distributions in $\mathcal M_1$ and $\I$ be a collection $\I$ of intervals in $[0,1]$.
	We say that (a) $\mathcal M$ is \emph{closed under
		concentration within $\I$} if  $F^{\I} \in  \mathcal M$ for all $F\in\M$;  (b) $\mathcal M$ is \emph{closed under
		concentration for all intervals} if for all $F\in \mathcal M$, we have $F^C\in \mathcal M$ for all intervals $C\subset [0,1]$; (c) $\mathcal M$ is \emph{closed under
		conditional expectation} if for all $F_X\in \mathcal M$, the distribution of any conditional expectation of $X$ is in $\mathcal M$.  
\end{definition}
 
 The relationship between the three properties of closedness in Definition \ref{def:con} is discussed in Propositions \ref{lem:lem-2} and \ref{lem:lem-4} below.
 Generally,  (c)$\Rightarrow$(b)$\Rightarrow$(a) if $\I$ is finite. 
Our main equivalence result is summarized in the following theorem. 
	\begin{theorem}\label{th:1}
	For $\mathcal M\subset \mathcal M_1$ and $h\in \mathcal H$, the following hold. 
\begin{enumerate}[(i)]
	\item	If  $h=\hat{h}$, i.e., $h$ is upper semicontinuous on $(0,1)$, and $\mathcal M $  is closed under
		concentration within $\I_h$,
		then 
		\begin{equation}\label{eq:eq-main}
			\sup_{F_Y\in \mathcal M } \rho_h(Y) = \sup_{F_Y\in \mathcal M}  \rho_{h^*}(Y).
		\end{equation}
	\item If $\mathcal M$  is closed under 
		concentration for all intervals, then \eqref{eq:eq-main} holds.
	\item	If $h=\hat{h}$, $\mathcal M$ is closed under concentration within $\I_h$, and the second supremum in \eqref{eq:eq-main} is attained by some $F\in \mathcal M$, then  $F^{\mathcal {I}_h}$  attains both suprema.
		\end{enumerate}
	\end{theorem}
	

%
	
Both suprema in \eqref{eq:eq-main} may be infinite, and this is discussed in Remark \ref{rem:infinite} in Appendix 
\ref{app:a5}.
	The proof of Theorem \ref{th:1} is   more technical  than   similar results in the literature because of the challenges arising from non-monotonicity, non-positivity, and discontinuity of $h$; see Figure \ref{fig:hat-h} for a sample of possible complications. In (ii),   $h$ does not need to be upper semicontinuous on $(0,1)$ for \eqref{eq:eq-main} to hold because closedness  under concentration for all intervals in (ii) is stronger than the condition in (i).

\begin{remark}
For $\mathcal{M}\subset\mathcal{M}_1$ and $h\in\H$, if $h=\hat{h}$ and $F^C\in\M$ for all $F\in\M$ and $C\in\I_h$, then the equivalence relation \eqref{eq:eq-main} also holds.   If $\I_h$ is finite, then this condition is  generally stronger than closedness under concentration within $\I_h$ in (i). 
\end{remark}

A natural question from Theorem \ref{th:1}  is whether our key condition of closedness under concentration  is necessary in some sense  for the equivalence  \eqref{eq:eq-main} to hold.\footnote{We thank an anonymous  referee for raising this question.} 
It is immediate to notice that adding any distributions $F_Z$ satisfying
$\rho_{h^*}(Z)<\sup_{F_Y\in\M}\rho_{h^*}(Y)$ 
 to the set $\M$ does not affect the equivalence,
 and therefore we turn our attention to the set of maximizers instead of the whole set $\M$.
In the next result, we show that closedness under concentration  within $\I_h$ of the set of maximizers of \eqref{eq:yproblem} is  necessary   for the equivalence \eqref{eq:eq-main} to hold.
\begin{theorem}\label{th:nec}
For $\mathcal{M}\subset\M_1$ and $h\in\H$ such that  $h\neq h^*$,  suppose that  the set   $\M_{\mathrm{opt}} $ of all maximizers of $\max_{F_Y\in\M}\rho_{h}(Y)$  is non-empty. 
If the equivalence \eqref{eq:eq-main} holds, i.e., 
$\sup_{F_Y\in \mathcal M } \rho_h(Y) = \sup_{F_Y\in \mathcal M}  \rho_{h^*}(Y)$, then $\M_{\mathrm{opt}}$ is closed under concentration within $\I_h$.
\end{theorem} 
 
If the equivalence \eqref{eq:eq-main} holds, then  
each $F\in\M_{\mathrm{opt}}$ also maximizes the problem $\sup_{F_Y\in \mathcal M}  \rho_{h^*}(Y)$.  
Conversely, if $h=\hat h$, then this condition and closedness of $\M_{\mathrm{opt}}$ under concentration within $\I_h$  together are necessary (by Theorem \ref{th:nec}) and sufficient (by Theorem \ref{th:1}) for the equivalence \eqref{eq:eq-main} to hold.   
If the maximizer $F$ of the original problem \eqref{eq:yproblem} is unique, then by Theorem \ref{th:nec}, $F$ must be equal to $F^{\I_h}$. 
The equivalence   \eqref{eq:eq-main} does not imply  closedness under concentration within $\I_h$ of the uncertainty set $\mathcal{M}$ itself; an example showing this is discussed in Remark \ref{rem:boundary}.  

\subsection{Some examples of distortion riskmetrics}
\label{sec:riskmetrics}

\begin{figure}[t]
		\begin{center}
			\includegraphics[width=0.4\textwidth]{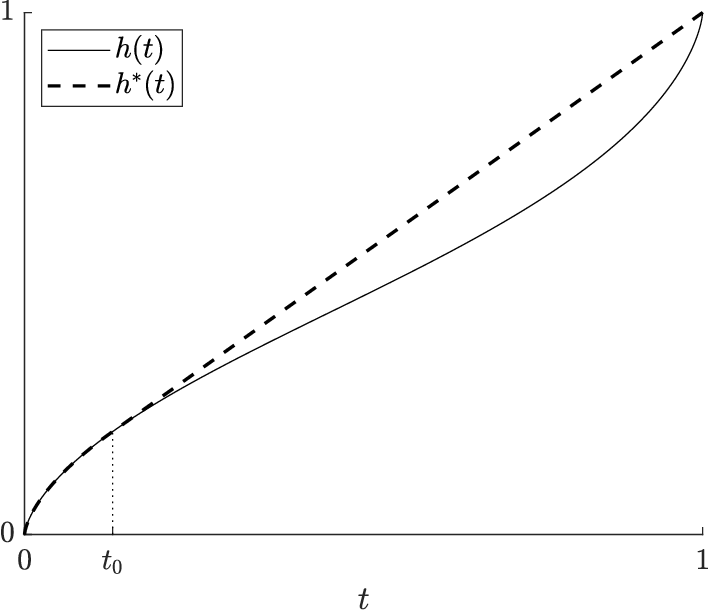} ~~~~
		   \includegraphics[width = 0.4\textwidth]{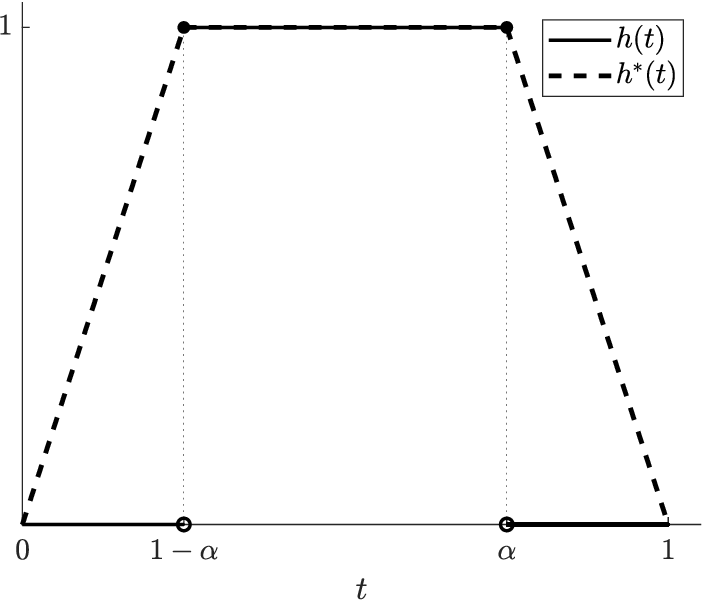}
		\end{center}
		\caption{Left panel: $h$  and $h^*$ for the TK distortion riskmetric   with $\gamma=0.7$ in Example \ref{ex:inverseS}; right panel: $h$ and $h^*$ for the inter-quantile range in Example \ref{ex:iqr}}
		\label{fig:TKdistortion}
		\end{figure}

\begin{figure}[htbp]
		    \centering
		   \includegraphics[width = 0.4\textwidth]{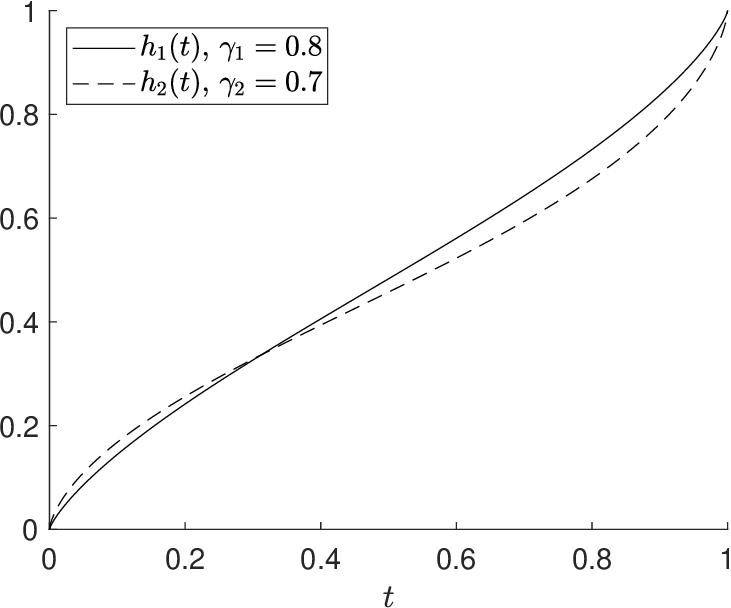} ~~~~
		   \includegraphics[width = 0.4\textwidth]{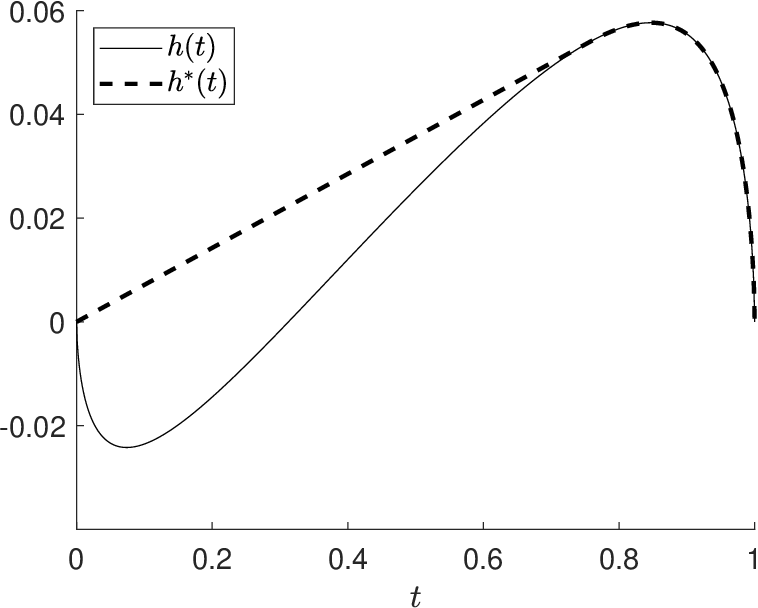}
		    \caption{Left panel: inverse-S-shaped distortion functions $h_1$ and $h_2$ in Example    \ref{ex:diff}; right panel: $h=h_1-h_2$ and $h^*$ of the same example}
		    \label{fig:distortion}
	\end{figure}

We provide a few examples of distortion riskmetrics $\rho_h$ commonly used in decision theory and finance, and obtain their corresponding set $\I_h$.
The Value-at-Risk (VaR) and the Expected Shortfall (ES) are the most popular risk measures in practice. We introduce them first,  followed by an inverse-S-shaped distortion function of \cite{TK92}.
	\begin{example}[VaR and ES]\label{ex:vares} For   $Y\in \L^0$, 
	  using the sign convention of \cite{MFE15},  VaR is defined as the left-quantile, and upper VaR ($\mathrm{VaR}^+$) is defined as the right-quantile; that is,
	$$\VaR_\alpha (Y)=F^{-1}_Y(\alpha),~~\alpha\in (0,1]~~~ \mbox{and} ~~~\VaR^+_\alpha (Y)=F^{-1+}_Y(\alpha),~~\alpha\in [0,1).$$   
	ES at level $\alpha$ is defined as 
	$$
	\ES_\alpha(Y)= \frac{1}{1-\alpha} \int_\alpha^1 \VaR_t(Y)\d t  ,~~\alpha\in(0,1),~Y\in \mathcal L^1.
	$$
	Both $\VaR_{\alpha}$ and $\ES_\alpha $ belong to the class of distortion riskmetrics.  
		Take $\alpha \in (0,1)$. Let $h (t) = \id_{(1-\alpha,1]} (t) $, $t\in[0,1]$. It follows that $h\in\H$ and $\hat{h}(t)=\id_{[1-\alpha,1]} (t)$, $t\in [0,1]$.  In this case, 
		$\rho_h =\VaR_\alpha $.
		Moreover, $h^*(t)=  \frac{t}{1-\alpha}\wedge 1$, $t\in [0,1]$ and $\rho_{h^*} =\ES_\alpha $.  Since $h^*$ and $\hat{h}$ differ on $(0,1-\alpha),$ we have $\I_h=\{(\alpha,1)\}$.
	\end{example} 
	\begin{example}[TK distortion riskmetrics]
		\label{ex:inverseS}
		The following function $h$ is an inverse-S-shaped distortion function (see also Figure \ref{fig:distortion}):
		\begin{equation}\label{eq:TKdistortion} h(t)=\frac{t^{\gamma}}{\left(t^{\gamma}+(1-t)^{\gamma}\right)^{1 / \gamma}},~~t\in[0,1],~\gamma\in(0,1).\end{equation}
		 Distortion riskmetrics with   distortion function  \eqref{eq:TKdistortion} are commonly used in behavioural economics and finance; see e.g., \cite{TK92}. 
		For simplicity, we call such distortion riskmetrics \emph{TK distortion riskmetrics}. Typical values of $\gamma$ are in $[0.5,0.9]$; see \cite{WG96}. For $h$ in \eqref{eq:TKdistortion}, it is clear that $h=\hat{h}$ on $[0,1]$ by continuity of $h$. We have $h^*\neq h$ on $(t_0,1)$, for some $t_0\in(0,1)$, and $h^*$ is linear on $[t_0,1]$. 
		Thus,  $\mathcal{I}_h=\{(0,1-t_0)\}$. An example of $h$ in \eqref{eq:TKdistortion} and its concave envelope $h^*$ are plotted in Figure \ref{fig:TKdistortion} (left).
	\end{example}

	For $h_1,h_2\in\H$, we write $h=h_1-h_2\in\H$ and consider the difference between two distortion riskmetrics, that is
		\begin{equation}
		\label{eq:diff}
			\rho_{h}=\rho_{h_1}-\rho_{h_2}.
		\end{equation}
		Such type of distortion riskmetrics measure the difference or disagreement between two utilities, risk attitudes, or capital requirements. Determining the upper and lower bounds, or the largest absolute values of such measures of disagreement, is of interest in practice but rarely studied in the literature. Note that $h_1-h_2$ is in general not monotone or concave even when $h_1$ and $h_2$ themselves have the specified  properties. Below we show some examples of distortion riskmetrics taking the form of \eqref{eq:diff}.
	\begin{example}[Inter-quantile range and inter-ES range]
	\label{ex:iqr}
		For $\alpha\in[1/2,1)$, we take $h_1(t)=\id_{[1-\alpha,1]} (t)$ and $h_2(t)=\id_{(\alpha,1]} (t) $, $t\in[0,1]$. It follows that $h(t)=h_1(t)-h_2(t)=\id_{\{1-\alpha\le t\le \alpha\}}$, $t\in[0,1]$, $\hat{h}=h$, and $$\rho_h(X)=F^{-1+}_X(\alpha)-F^{-1}_X(1-\alpha),~~X\in\L^0.$$
		Correspondingly, we have $h^*(t)=t/(1-\alpha)\wedge 1+(\alpha-t)/(1-\alpha)\wedge 0$, $t\in[0,1]$, and
		$$\rho_{h^*}(X)=\ES_\alpha(X)+\ES_\alpha(-X),~~X\in\L^1.$$
		This distortion riskmetric $\rho_h$ is called an inter-quantile range and $\rho_{h^*}$ is called an inter-ES range. As the distortion functions $h^*$ and $\hat{h}$ differ on the open intervals $(0,1-\alpha)$ and $(\alpha,1)$, we have $\mathcal{I}_h=\{(\alpha,1),(0,1-\alpha)\}$. The distortion functions $h$ and $h^*$ are displayed in Figure \ref{fig:TKdistortion} (right).
	\end{example}	
	
	\begin{example}[Difference of two inverse-S-shaped distortion functions]
	\label{ex:diff}
		We take $h_1$ and $h_2$ to be the inverse-S-shaped distortion functions in \eqref{eq:TKdistortion}, with parameters $\gamma_1=0.8$ and $\gamma_2=0.7$, respectively. By calculation, the function $h=h_1-h_2$ is convex on $[0,0.3770]$, concave on $[0.3770,1]$, and as seen in Figure \ref{fig:distortion} not monotone. The concave envelope $h^*$ is linear on $[0,0.7578]$ and $h^*=h$ on $[0.7578,1]$. Thus, we have $\mathcal{I}_h=\{(0.2422,1)\}$. The graphs of the distortion functions $h_1$, $h_2$, $h$, and $h^*$ are displayed in Figure \ref{fig:distortion}.
	\end{example}
	
	The functions in $\mathcal H$ are a.e.~differentiable, and for an absolutely continuous function $h\in \mathcal H$, let $h'$ be a (representative) function on $[0,1]$ that is a.e.~equal to the derivative of $h$.  
	If $h\in \mathcal H$ is left-continuous or $\VaR_t(Y)$ is continuous with respect to $t\in(0,1)$, the risk measure $\rho_h$ in \eqref{eq:disrep} has representation 
	\begin{equation}\label{eq:disrep2} \rho_h(Y) =\int_0^1 \VaR_{1-t}(Y) \d h(t),~~Y\in \L^p;\end{equation} 
	see Lemma 1 of \cite{WWW20}. If $h\in \mathcal H$ is absolutely continuous it holds
	\begin{equation}\label{eq:disrep3} 
		\rho_{h}(Y) =\int_0^1 \VaR_{1-t}(Y) h'(t)\d t, ~~Y\in \L^p.
	\end{equation}  
	Another example of a recently introduced distortion riskmetric with concave distortion function  may be of independent interest in risk management.
	\begin{example}[Second-order superquantile]
	\label{ex:ssq}
		As introduced by \cite{RR18},  a second-order superquantile  is defined as
		$$\mathrm{SSQ}_\alpha(Y)=\frac{1}{1-\alpha}\int^1_\alpha\ES_t(Y)\d t,~~\alpha\in(0,1),~Y\in\L^2.$$  
		By Theorem 2.4 of \cite{RR18}, $\mathrm{SSQ}_\alpha $ is a  distortion riskmetric  with a concave distortion function $h$ given by
		$$h(t)=\left\{\begin{array}{ll}
		\frac{t}{1-\alpha}\left(1+\log\frac{1-\alpha}{t}\right),& 0\le t< 1-\alpha,\\
		1 , & 1-\alpha\le t\le 1.
		\end{array}\right.$$ 	
Clearly, $\mathrm{SSQ}_ \alpha \ge \ES_\alpha$. The difference $\mathrm{SSQ}_ \alpha-\ES_\alpha$ between second-order superquantile and ES, which has a similar interpretation as $\ES_\alpha-\VaR_\alpha$, is a distortion riskmetric with a non-concave and non-monotone distortion function $g$, and the set $\mathcal I_{g}$ contains a single interval of the form $(0,\beta)$ for some $\beta\in [\alpha,1)$.
	\end{example}

\subsection{Closedness under concentration for all intervals}

In this section, we present some technical results and specific examples about closedness under concentration for all intervals and under conditional expectation. 
	The proposition below clarifies  the relationship between closedness under concentration for all intervals 
	and closedness under
	conditional expectation.

	\begin{proposition}\label{lem:lem-2}
		Closedness under
		conditional expectation implies closedness under concentration for all intervals, but the converse is not true.
	\end{proposition}


	
	\begin{example} \label{ex:main} 
 We present 6   classes of sets $\mathcal M$ that are closed under conditional expectation, and hence also under concentration for all intervals.  
		\begin{enumerate} 
			\item (Moment conditions) For $p>1$, $m\in \R$, and $v>0$,   the set $$
			\mathcal M(p,m,v) = \{ F_Y\in \mathcal M_p: \E[Y]=m, ~\E[|Y-m|^{p}]\le v^{p}\}
			$$ 
			is closed under conditional expectation by Jensen's inequality.
			The set $ \mathcal M(p,m,v)$ corresponds to distributional uncertainty with moment information, and the setting $p=2$  (mean and   variance constraints) is the most commonly studied.
			
			\item (Mean-covariance conditions) For $n\in\N$, $\mathbf{a}\in\R^n$, $\bm{\mu}\in\R^n$, and $\Sigma\in \R^{n\times n}$ positive semidefinite, let
			$$\M^{\mathrm{mv}}(\mathbf{a},\bm{\mu},\Sigma)=\{F_{\mathbf{a}^{\top}\mathbf{X}}\in\M_2:F_{\mathbf{X}}\in\M^n_2,~\E[\mathbf{X}]=\bm{\mu},~\mathrm{var}(\mathbf{X})\preceq\Sigma\},$$  
			where $\mathbf{X} =(X_1,\dots,X_n)$, $\E[\mathbf{X}] =(\E[X_1],\dots,\E[X_n])$,   $\mathrm{var}(\mathbf{X})$ is the covariance matrix of $\mathbf{X}$,
			and  $B^\prime\preceq B$ means that the matrix $B-B^\prime$ is positive semidefinite for two positive semidefinite symmetric matrices $B$ and $B^\prime$.
			With a simple verification in   Appendix \ref{app:main}, $\M^{\mathrm{mv}}(\mathbf{a},\bm{\mu},\Sigma)=\M(2,\mathbf{a}^{\top}\bm{\mu},(\mathbf{a}^{\top}\Sigma\mathbf{a})^{1/2})$.
			
			\item (Convex function conditions)  For $n\in\N$, $\mathbf{a}\in \R^n$, $K\subset\N$,  a collection $\mathbf{f}=(f_k)_{k\in  K}$ of convex functions  on $\R^n$,  and a vector $\mathbf{x}=(x_k)_{k\in  K}\in \R^{|K|}$, let  $$
			\mathcal M^\mathbf{f} (\mathbf{a},\mathbf{x}) = \{ F_{\mathbf{a}^{\top}\mathbf{X}}\in \mathcal M_1: \E[f_k(\mathbf{X})]\le x_k~\mbox{for all} ~k\in K\}. 
			$$  
			The set $ \mathcal M ^\mathbf{f}$ corresponds to distributional uncertainty with constraints on expected losses or test functions. The set $\mathcal M^{\mathbf {f}}$  includes $\mathcal M(p,m,v)$ as a special case.
			
			\item (Distortion  conditions) 
			For $K\subset\N$, a collection   $\mathbf{h}=(h_k)_{k\in  K}\in   (\mathcal H^*) ^{|K|} $ and a vector $\mathbf{x}=(x_k)_{k\in  K}\in \R^{|K|} $,  let $$
			\mathcal M^\mathbf{h}(\mathbf{x})  = \{ F_Y\in \mathcal M_1:   \rho_{h_k}(Y)\le x_k~\mbox{for all} ~k\in K\}. 
			$$  
			The set $ \mathcal M ^\mathbf{h}$ corresponds to distributional uncertainty with constraints on preferences modeled by convex dual utilities. 
			\item (Convex order conditions) 
			For $K\subset\N$ and a collection of random variables $\mathbf{Z}=(Z_k)_{k\in K}\in (\L^1)^{|K|}$, let $$
			\mathcal M^{\rm cx}(\mathbf{Z})  = \{ F_Y\in \mathcal M_1 :  Y\le_{\rm cx} Z_k~\mbox{for all} ~k\in K\},
			$$
			  where $\le_{\rm cx}$ is the inequality in convex order.\footnote{Precisely, we write $G\le_{\rm cx} (\le_{\rm icx}) \,F$ if  $\int \phi \d G \le \int \phi  \d F$ for all (increasing) convex functions $\phi$ such that the above two integrals are well defined. }  
			Similar to the above two examples,  $ \mathcal M ^{\rm cx} (\mathbf{Z})$ is closed under conditional expectation (cf.~Remark \ref{rem:cxorder} in Appendix \ref{app:a5}).
			\item (Marginal conditions) 
			For given univariate distributions $F_1,\dots,F_n\in \mathcal M_1$, 
			let $$
			\mathcal M^S(F_1,\dots,F_n)  = \{ F_{X_1+\dots+X_n} \in \mathcal M_1 : X_i \sim F_i,~i=1,\dots,n\}. 
			$$  
			In other words, $\mathcal M^S$ is the set of all possible aggregate risks $X_1+\dots+X_n$ with given marginal distributions of $X_1,\dots,X_n$; see \cite{EWW15} for some results  on $\mathcal M^S$.
			Generally, $\mathcal M^S$ is not closed under concentration for all intervals or conditional expectation, since closedness under concentration for all intervals is stronger than joint mixability \citep{WW16}. 
			In the special case where $F_1=\dots=F_n=\mathrm{U}[0,1]$,    
			Proposition 1 and Theorem 5 of \cite{MWW19} imply that $\mathcal M^S $ is   closed under conditional expectation if and only if $n\ge 3$.
		\end{enumerate}
	\end{example}

	\begin{remark}\label{rem:boundary}
		The uncertainty set $\mathcal{M}(p,m,v)$ of the moment condition in Example \ref{ex:main} can be restricted to the set
		$$\overline{\mathcal {M}}(p,m,v)=\{ F_Y\in \mathcal M_p: \E[Y]=m, ~\E[|Y-m|^{p}]= v^{p}\},$$
		which is the ``boundary" of $\mathcal{M}(p,m,v)$. For    $\M=\mathcal{M}(p,m,v)$, the suprema   on both sides of \eqref{eq:eq-main} are obtained by some distributions in $\overline{\mathcal {M}}(p,m,v)$; see Theorem \ref{th:th1}. As a direct consequence, we get   $$\sup_{F_Y\in\overline{\mathcal {M}}(p,m,v)}\rho_{h^*}(Y)=\sup_{F_Y\in{\mathcal {M}}(p,m,v)}\rho_{h^*}(Y)=\sup_{F_Y\in{\mathcal {M}}(p,m,v)}\rho_h(Y)=\sup_{F_Y\in\overline{\mathcal {M}}(p,m,v)}\rho_h(Y).$$ Hence, equivalence holds  even though $\overline{\mathcal {M}}(p,m,v)$ is not closed under concentration for any interval. 
		By Theorem \ref{th:nec}, the set of optimizers is closed under concentration within $\I_h$ for each $h\in \mathcal H$. 
	\end{remark}

	For a distribution $F\in\M_1$ and a collection $\I$ of disjoint intervals in $[0,1]$, we have the following result regarding to the distribution $F^\mathcal{I}$.
	\begin{proposition}\label{lem:lem-4}
		Let $\I$ be a collection of disjoint intervals in $[0,1]$ and $\mathcal M$ be a set of distributions. If $\M$ is closed under concentration for all intervals and $\I$ is finite, or $\mathcal M$ is closed under conditional expectation, then  $\M$ is closed under concentration within $\I$.
	\end{proposition}

If $\I$ is infinite, closedness under concentration for all intervals may not be sufficient for closedness under concentration within  $\I$; 
see Remark \ref{rem:counterex} in Appendix \ref{app:a5} for a technical explanation.
An infinite $\I_h$ does not appear for any distortion riskmetrics in practice.


	\subsection{Examples of closedness under concentration within $\I$ but not for all intervals}
	 \label{sec:example}

	In practice, it is much easier to check closedness under concentration within a specific collection of intervals $\I$ than closedness under concentration for all intervals or under conditional expectation. In this section, we show several examples for closedness under concentration within some $\I$.
	
	For 
	distortion functions $h$ such that $\I_h=\{(p,1)\}$ (resp.~$\I_h=\{(0,p)\}$) for some $p\in (0,1)$,
	the result in Theorem \ref{th:1} {(i)} only requires $\mathcal M$ to be 
	closed under concentration within $\{(p,1)\}$ (resp.~$\{(0,p)\}$).
	Such distortion functions  
	include the
	inverse-S-shaped distortion functions in \eqref{eq:TKdistortion},  those of $\VaR_p$, and $\VaR^+_p$, and that of the difference between the second-order superquantile and ES in Example \ref{ex:ssq}. 
	Below we present some more concrete examples.
	

	\begin{example}[$\M$ has two elements]
	Let $p\in (0,1)$ and  $\mathcal M=\{\mathrm{U}[0,1],\, p\delta_{p/2}+ (1-p)\mathrm{U}[p,1]\}$ where $\delta_{p/2}$ is the point-mass at $p/2$. We can check that $\mathcal M$ is closed under concentration within $\{(0,p)\}$ but $\mathcal M$ is not closed under concentration for all intervals. Indeed, any set closed under concentration for all intervals and containing	$\mathrm{U}[0,1]$ has infinitely many elements. 
 In general, a finite set which contains any non-degenerate distribution is not closed under conditional expectation in an atomless probability space, since there are infinitely many possible distributions for the conditional expectation of a given non-constant random variable. 	 Another similar example that is closed under concentration within $\{(0,p)\}$ is the set of all possible distributions of the sum of several Pareto risks; see Example 5.1 of \cite{WXZ19}.
	\end{example}
	\begin{example}[VaR and ES]
		As we see from Example \ref{ex:vares},
		if
		$\rho_h=\VaR^+_\alpha$ for some $\alpha \in (0,1)$,
		then $\rho_{h^*}$ is $\ES_\alpha$ and  $\I_h=\{(\alpha,1)\}$.
		Theorem \ref{th:1} {(i)} implies that if $\mathcal M$ is closed under concentration within $\{(\alpha,1)\}$, then  $$
		\sup_{F_Y\in \mathcal M } \VaR^+_\alpha(Y) = \sup_{F_Y\in \mathcal M}  \ES_\alpha(Y).
		$$ 
		This observation leads to (with some modifications) the main results in \cite{WBT15} and \cite{LSWY18} on the equivalence between VaR and ES.
	\end{example}
	\begin{example}[TK distortion riskmetric]
	\label{ex:9}
		If we take $h$ to be an inverse-S-shaped distortion function in \eqref{eq:TKdistortion}, then $\mathcal{I}_h=\{(0,1-t_0)\}$ for some $t_0\in(0,1)$, and $\rho_h$ is the TK distortion riskmetric. As a direct consequence of Theorem \ref{th:1} {(i)}, if $\M$ is closed under concentration within $\{(0,1-t_0)\}$, then
		$$\sup_{F_Y\in\M}\rho_h(Y)=\sup_{F_Y\in\M}\rho_{h^*}(Y).$$
		This result implies Theorem 4.11 of \cite{WXZ19} on the robust risk aggregation problem based on dual utilities with inverse-S-shaped distortion functions.
	\end{example}

\begin{example}[Wasserstein ball, $1$-dimensional]
\label{ex:wassball}
Optimization problems under the uncertainty set of a Wasserstein ball are common in literature when quantifying the discrepancy between a benchmark distribution and alternative scenarios; see e.g., \cite{BM19}. We discuss the application of the concept of concentration to optimization with Wasserstein distances. For $p\ge 1$ and $F,G\in\M_p$, the \emph{$p$-Wasserstein distance} between $F$ and $G$ is defined as
$$W_p(F,G)=\left(\int^1_0\left|F^{-1}(u)-G^{-1}(u)\right|^p\d u\right)^{1/p}.$$
For $\epsilon\ge 0$, the uncertainty set of an \emph{$\epsilon$-Wasserstein ball} around a benchmark distribution $\widetilde{G}\in\M_p$ is given by
$$\mathcal{M}(\widetilde{G},\epsilon)=\{F\in\M_p:W_p(F,\widetilde{G})\le\epsilon\}.$$
Suppose that the benchmark distribution $\widetilde{G}$ has a quantile function that is constant on each element in some collection of disjoint intervals $\widetilde{\I}\subset[0,1]$. As shown in Appendix \ref{app:main}, $\M(\widetilde{G},\epsilon)$ is closed under concentration within $\I$ for all $\I\subset\widetilde{\I}$. 
Using this  closedness property  and  Theorem \ref{th:1} (i),
  the equivalence 
\begin{equation}\label{eq:wassball} \sup_{F_Y\in\M(\widetilde{G},\epsilon)}\rho_h(Y)=\sup_{F_Y\in\M(\widetilde{G},\epsilon)}\rho_{h^*}(Y) \end{equation}
holds for all $h\in\H$ such that $\I_h\subset\widetilde{\I}$.
\end{example}

\begin{remark}
In general, if the quantile function $\widetilde{G}$ in Example \ref{ex:wassball} is not constant on some interval in $\widetilde{\mathcal I}$, then $\M(\widetilde{G},\epsilon)$ is not necessarily closed under concentration within $\widetilde{\mathcal I}$, and the   equivalence \eqref{eq:wassball} may not hold. For instance, the worst-case $\VaR_\alpha$ over $\M(\widetilde{G},\epsilon)$ is generally different from the worst-case $\ES_\alpha$  over $\M(\widetilde{G},\epsilon)$ as obtained in Proposition 4 of \cite{LMWW22}.  We also refer to \cite{BPV20} who consider a Wasserstein ball together with moment constraints.
\end{remark}

\begin{example}[Wasserstein ball, $n$-dimensional]\label{ex:11}
For $n\in\N$, $p\ge 1$, $a\ge 1$ and $F,G\in\M^n_p$, the $p$-Wasserstein distance on $\R^n$ between $F$ and $G$ is defined as
$$W^n_{a,p}(F,G)=\inf_{\mathbf{X}\sim F,~\mathbf{Y}\sim G}(\E[\Vert \mathbf{X}-\mathbf{Y}\Vert ^p_a])^{1/p},$$
where $\Vert \cdot\Vert _a$ is the $\mathcal L^a$-norm on $\R^n$. Similarly to the $1$-dimensional case, for $\epsilon\ge 0$, an $\epsilon$-Wasserstein ball on $\R^n$ around a benchmark distribution $\widetilde{G}\in\M^n_p$ is defined as
$$\mathcal{M}^n(\widetilde{G},\epsilon)=\{F\in\M^n_p:W^n_{a,p}(F,\widetilde{G})\le \epsilon\}.$$
In a portfolio selection problem, we consider the worst-case riskmetric of a linear combination of random losses. For $\epsilon\ge 0$, $\mathbf{w}\in[0,\infty)^n$, $p>1$, $a>1$ and $\mathbf{Z}\in(\L^p)^n$, as shown in Appendix \ref{app:main}, the uncertainty set
$$\{F_{\mathbf{w}^{\top}\mathbf{X}}\in\M_p:F_\mathbf{X}\in\mathcal{M}^n(F_{\mathbf{Z}},\epsilon)\}$$
is closed under concentration within $\{(0,t)\}$ for all $t\le p_0$. 
For a practical example, assume that an investor holds a portfolio of bonds (for simplicity, assume that they have the same maturity). The loss vector  $\mathbf X\ge\mathbf 0$  from this portfolio at maturity has an estimated benchmark loss distribution $\widetilde{G}$, and  the probability of no default  from these bonds  (i.e., $\mathbf X=\mathbf 0$)  is estimated as $p_0>0$ (usually quite large). Suppose that the investor uses a distortion riskmetric with an inverse-S-shaped distortion function $h$ given in \eqref{eq:TKdistortion} of Example \ref{ex:inverseS}, and considers a Wasserstein ball around $\widetilde G$ with radius $\epsilon$.
Note that $\I_h=\{(0,t)\}$ for some $t\in(0,1)$ from Example \ref{ex:9}. 
By Theorem \ref{th:1} (i), we obtain an equivalence result on the worst-case riskmetrics for the portfolio with weight vector $\mathbf w$,
$$\sup_{F_\mathbf{X}\in\M^n(\widetilde{G},\epsilon)}\rho_h(\mathbf{w}^{\top}\mathbf{X})=\sup_{F_\mathbf{X}\in\M^n(\widetilde{G},\epsilon)}\rho_{h^*}(\mathbf{w}^{\top}\mathbf{X}),$$
whenever $t\in(0,p_0]$.
\end{example}

\begin{example}[Optimal hedging strategy]
\label{ex:opt_hedge}

Suppose that an investor is willing to hedge her random loss $X$ only when it exceeds some certain level $l\in\R$. Mathematically, for a fixed $X\in\L^1$ continuously distributed on $(F^{-1}_X(p_0),F^{-1}_X(1))$ such that $\P(X\le l)=p_0$ for some $p_0\in(0,1)$ and $l\in\R$, define the set of measurable functions $$\mathcal{V}=\{V:\R\to\R\mid x\mapsto x-V(x) \text{ is increasing},~  V(x)=0\text{ for all }x\le l\}$$ representing possible hedging strategies. Let $g:\R\to\R$ be an increasing and convex function. The final payoff obtained by a hedging strategy $V\in\mathcal{V}$ is given by $X-V(X)+g(\E[V(X)])$, where $g(\E[V(X)])$ is a fixed cost of the hedging strategy that depends on the expected value of $V(X)$ calculated by a risk-neutral seller in the market using the same probability measure $\P$.
As shown in Appendix \ref{app:main}, the action set in this optimization problem,
$$\mathcal{M}=\{F_{X-V(X)+g(\E[V(X)])}\in\M_1:V\in\mathcal{V}\},$$
is closed under concentration within $\{(p,1)\}$ for all $p\in[p_0,1)$.
On the other hand, it is obvious that $\mathcal{M}$ is not closed under concentration for all intervals or closed under conditional expectation since the quantiles of the distributions in $\mathcal{M}$ are fixed beyond the interval $(p_0,1)$. 
The above closedness under concentration property allows us to use 
Theorem \ref{th:1}  to convert the optimal hedging problem for $\rho_h$ with 
an inverse-S-shaped distortion function $h$ as in \eqref{eq:TKdistortion} to a convex version $\rho_{h^*}$.


\end{example}

\begin{example}[Risk choice]

Suppose that an investor is faced with a random loss $X\in\L^1$. The distortion function $h$ of her riskmetric is inverse-S-shaped with $\mathcal{I}_{-h}=\{(p,1)\}$ for some $p\in(0,1)$. Suppose that $p$ is known to the seller. Since the investor is averse to risk for large losses, the seller may provide her with the option to stick to the initial investment or to convert the upper part of the random loss into a fixed payment to avoid large loss. Specifically, we consider the   set $\mathcal{M}=\{F_X,F_X^{(p,1)}\}$ containing two elements, where $\P(X\le u)=p$ for some $u\in\R$. It is clear that $\mathcal{M}$ is closed under concentration within $\{(p,1)\}$ but not closed under conditional expectation. 
We assume that the costs of the two investment strategies are calculated by expectation and thus are the same. By (i) of Theorem \ref{th:1}, it  follows that the risk minimization problem satisfies
$$\min_{F_Y\in\mathcal{M}}\rho_h(Y)=\min_{F_Y\in\mathcal{M}}\rho_{h_*}(Y)=\rho_{h_*}(X),$$
where the last equality follows from Theorem 3 of \cite{WWW20}. By (iii) of Theorem \ref{th:1}, we further have the minimum of the original problem $\min_{F_Y\in\mathcal{M}}\rho_h(Y)$ is obtained by $F_X^{(p,1)}$; intuitively, the investor will choose
to convert the upper part of her loss into a fixed payment. 


%

\end{example}

\subsection{Atomic probability space}

The definition of closedness under concentration  in Definition \ref{def:con} requires the assumption of an atomless probability space since a uniform random variable is used in the setup. It may be of practical interest in some economic and optimization settings to assume a finite probability space. In this section, we let the sample space be $\Omega_n=\{\omega_1,\dots,\omega_n\}$ for $n\in\N$ and the probability measure $\P_n$ be such that $\P_n(\omega_i)=1/n$ for all $i=1,\dots,n$ (such a space is called \emph{adequate} in economics).
The possible distributions in such a probability space are supported by at most $n$ points each with probability a multiple of $1/n$, and we denote by $\mathcal M_{[n]}$ the set of these distributions.

Define the collection of intervals $\mathcal{I}_n=\{(j/n,k/n]:j,k\in\N\cup\{0\},~j<k\le n\}$. We say a set of  distributions $\mathcal{M} \subset \M_{[n]}$   is \emph{closed under grid concentration} within $\mathcal{I}\subset\I_n$ if for all $F\in\mathcal{M}$, the distribution of the random variable
$$F^{-1}(U_n)\id_{\{U_n\notin \bigcup_{C\in\I} C\}}+\sum_{C\in\I}\E[F^{-1}(U_n)|U_n\in C]\id_{\{U_n\in C\}}$$ 
is also in $\mathcal{M}$, where $U_n$ is a random variable such that $U_n(\omega_i)=i/n$ for all $i=1,\dots,n$. 
For a distribution $F$ with finite mean and $(a,b]\in\I_n$, it is straightforward that
		the left-quantile function   of 
		$F^{(a,b]}$ is given by \eqref{eq:quantile}.
The following equivalence result holds with additional assumption $\I_h\subset \I_n$. The proof  can be obtained directly from that of Theorem \ref{th:1}.

\begin{proposition}\label{prop:discrete}
Let $\mathcal M \subset \M_{[n]}$  and $h\in \mathcal H$.
	If  $h=\hat{h}$, $\I_h\subset\I_n$ and $\mathcal M $  is closed under
		grid concentration within $\I_h$,
		then 
		\begin{equation*}
			\sup_{F_Y\in \mathcal M } \rho_h(Y) = \sup_{F_Y\in \mathcal M}  \rho_{h^*}(Y).
		\end{equation*}
\end{proposition}

We note that the condition $\I_h\subset \I_n$ in Proposition \ref{prop:discrete} is satisfied by all distortion functions $h$ which are linear (or constant) on each of $((j-1)/n,j/n]$, $j=1,\dots,n$.
It is common to assume such a distortion function $h$ in an adequate probability space of $n$ states, since any distribution function can only take values in $\{j/n:j=0,\dots,n\}$.

\section{Multi-dimensional setting}
\label{sec:multi}

Our main equivalence results in Theorems \ref{th:1} and \ref{th:nec} are stated under the context of one-dimensional random variables.  In this section, we  discuss their generalization to multi-dimensional framework with a few additional steps. 

In the multi-dimensional setting, closedness under concentration is not easy to define,  as  quantile functions are not naturally defined for multivariate distributions.
Nevertheless, closedness under conditional expectation can be analogously formulated. 
 {For $n\in\N$,} we say that  $\mathcal M\subset\mathcal M^n$ is \emph{closed under
		conditional expectation}, if for all $F_\mathbf{X}\in \mathcal M$, the distribution of any conditional expectation of $\mathbf{X}$ is in $\mathcal M$.  The following theorem states the multi-dimensional version of our main equivalence result using  closedness under conditional expectation.

\begin{theorem}\label{th:multi}
	For $\widetilde{\mathcal M}\subset \mathcal M^n_1$, increasing function $h\in \mathcal H$ and $f:A \times \R^n\to\R$ concave in the second argument, if $\widetilde{\mathcal M}$  is closed under 
		conditional expectation, then for all $\mathbf{a}\in A$,
		\begin{equation}\label{eq:multi}
			\sup_{F_\mathbf{X}\in \widetilde{\mathcal M} } \rho_h(f(\mathbf{a},\mathbf{X})) = \sup_{F_\mathbf{X}\in \widetilde{\mathcal M} } \rho_{h^*}(f(\mathbf{a},\mathbf{X})). 
		\end{equation}
		If $h=\hat{h}$ and the second supremum in \eqref{eq:multi} is attained by some $F_{\mathbf{X}}\in \widetilde{\mathcal {M}}$, then  $F^{\mathcal {I}_h}_{f(\mathbf{a},\mathbf{X})}$  attains both suprema.
	Moreover, if $f$ is linear in the second component, then \eqref{eq:multi} holds for all $h\in \mathcal H$ (not necessarily monotone). 
\end{theorem}

\begin{remark}
If we assume that $f$ is convex (instead of concave) in the second argument in  Theorem \ref{th:multi} and keep the other assumptions, then
for an increasing $h$,	$$\inf_{F_\mathbf{X}\in \widetilde{\mathcal M} } \rho_h(f(\mathbf{a},\mathbf{X})) = \inf_{F_\mathbf{X}\in \widetilde{\mathcal M} } \rho_{h_*}(f(\mathbf{a},\mathbf{X})).$$ 
	This statement follows by noting   $\rho_{-h}=-\rho_h$. The case of a decreasing $h$ is similar.
 \end{remark}
 
Theorem \ref{th:multi} is similar to Theorem 3.4 of \cite{CLM20} which states the equivalence \eqref{eq:multi} for increasing $h$ and a specific set $\widetilde \M$ which is a special case in Example \ref{ex:multi} below. 
In contrast, our result applies to non-monotone $h$ (with an extra condition on $f$), more general set $\widetilde \M$, and 
also  the infimum problem.
The setting of a function $f$ linear in the second argument
often appears in portfolio selection problems where $f(\mathbf a,\mathbf X)= \mathbf a^\top \mathbf X$; see Example \ref{ex:11} and Section \ref{sec:app}.

\begin{example}\label{ex:multi}
Similarly to Example \ref{ex:main}, we give examples of sets of multi-dimensional distributions closed under conditional expectation.
\begin{enumerate}
\item (Convex function conditions) For $n\in\N$, a convex set $B\subset\R^n$, set $\Psi$ of convex functions on $\R^n$,  and a mapping $\pi:\Psi\to\R$, let 
$$\widetilde{\M}(B,\Psi,f)=\{F_{\mathbf{X}}\in\M^n_1:\P(\mathbf{X}\in B)=1,~\E[\psi(\mathbf{X})]\le \pi(\psi)\text{ for all }\psi\in\Psi\}.$$
It is clear that $\widetilde{\M}(B,\Psi,f)$ is closed under conditional expectation due to Jensen's inequality. 
The uncertainty set proposed by \cite{DAY14} and used in Theorem 3.4 of \cite{CLM20}
can be obtained as a special case of this setting by taking $\Psi=\{f_1,\dots,f_n\}\cup\{g_1,\dots,g_n\}\cup\Phi$, where $f_i:(x_1,\dots,x_n)\mapsto x_i$, $g_i:(x_1,\dots,x_n)\mapsto -x_i$ for all $i=1,\dots,n$, and  $\Phi$ is a set of convex functions.
 The specification for $\pi$ is that $\pi(f_i)=m_i\in\R$, $\pi(g_i)=-m_i$, $\pi(\phi)=0$ for all $i=1,\dots,n$, $\phi\in\Phi$. 

\item (Distortion conditions) For $n\in\N$, $K\subset\N$, $\mathbf{a}=(\mathbf a_k)_{k\in K}\in \R^{n\times |K|}$,   $\mathbf{h}=(h_k)_{k\in  K}\in   (\mathcal H^*) ^{|K|} $ and  $\mathbf{x}=(x_k)_{k\in  K}\in \R^{|K|} $,  the set $$
			\widetilde{\mathcal M}^\mathbf{h}(\mathbf{a},\mathbf{x})  = \{ F_\mathbf{X}\in \mathcal M^n_1:   \rho_{h_k}(\mathbf{a}_k^{\top}\mathbf{X})\le x_k~\mbox{for all} ~k\in K\}
			$$  
is closed under conditional expectation. In portfolio optimization problems, this setting incorporates distributional uncertainty with constraints on convex distortion risk measures of the   total loss. In particular, optimization with the riskmetrics chosen as ES  is common in the literature; see e.g., \cite{RU02}, where  ES is called  CVaR.

\item (Convex order conditions) 
			For $n\in\N$ and random vectors $\mathbf{Z}_k\in(\mathcal{L}^1)^n$, $k\in K\subset \N$, we naturally extend from part 5 of Example \ref{ex:main} and obtain that the set $$
			\widetilde{\mathcal M}^{\rm cx}(\mathbf{Z})  = \{ F_\mathbf{X}\in \mathcal M^n_1 :  \mathbf{X}\le_{\rm cx} \mathbf{Z}_k\mbox{ for all }k\in K\}
			$$  
			is closed under conditional expectation. 
\end{enumerate}
\end{example}

Next, we discuss a multi-dimensional   problem setting involving concentrations of marginal distributions. For $n\in\N$, we assume that marginal distributions of an $n$-dimensional distribution in $\mathcal{M}^n_1$ are uncertain and are in some sets $\mathcal{F}_1,\dots,\mathcal{F}_n\subset\mathcal{M}_1$. For $F_1,\dots,F_n\in\mathcal{M}_1$, define the set
$$\mathcal{D}(F_1,\dots,F_n)=\{\text{cdf of }(X_1,\dots,X_n):X_i\sim F_i,~i=1,\dots,n\},$$
which is the set of all possible joint distributions with specified marginals; see \cite{EWW15}.
For $\mathbf{a}\in A$, $h\in\H$ and $\mathcal{F}_1,\dots,\mathcal{F}_n\subset\M_1$, the worst-case distortion riskmetric can be represented as
	\begin{equation}\label{eq:RRA}
	\sup_{F_{\mathbf{X}}\in\mathcal{D}(F_1,\dots,F_n)}\sup_{F_1\in\mathcal{F}_1,\dots,F_n\in\mathcal{F}_n}\rho_{h}(f(\mathbf{a},\mathbf{X})).
	\end{equation}
	The outer problem of \eqref{eq:RRA} is a robust risk aggregation problem (see \cite{EPR13, EWW15} and  {item 6} of Example \ref{ex:main}), which is typically nontrivial in general when $h$ is not concave. With additional uncertainty of the marginal distributions, \eqref{eq:RRA} can be converted to a convex problem given that $\mathcal{F}_1,\dots,\mathcal{F}_n$ are closed under concentration.
\begin{theorem}\label{prop:RRA}
	For $\mathcal{F}_1,\dots,\mathcal{F}_n\subset \mathcal M_1$, increasing $h\in \mathcal H$ with $h=\hat{h}$, and $f:A\times \R^n\to\R$ increasing, supermodular and positively homogeneous in the second argument, if  $\mathcal F_1,\dots,\mathcal{F}_n $  are closed under
		concentration within $\I_h$, then the following hold.\footnote{For a function $f:\R^n\to\R$, we say $f$ is \emph{supermodular} if $f(\mathbf{x})+f(\mathbf{y})\le f(\mathbf{x}\wedge \mathbf{y})+f(\mathbf{x}\vee\mathbf{y})$ for all $\mathbf{x},\mathbf{y}\in\R^n$; $f$ is \emph{positively homogeneous} if $f(\lambda\mathbf{x})=\lambda f(\mathbf{x})$ for all $\lambda\ge 0$ and $\mathbf{x}\in\R^n$.}
\begin{enumerate}[(i)]
	\item	For all $\mathbf{a}\in A$,
		\begin{equation}\label{eq:eq-multi}
			\sup_{F_{\mathbf{X}}\in\mathcal{D}(F_1,\dots,F_n)}\sup_{F_1\in\mathcal{F}_1,\dots,F_n\in\mathcal{F}_n}\rho_{h}(f(\mathbf{a},\mathbf{X}))=\sup_{F_{\mathbf{X}}\in\mathcal{D}(F_1,\dots,F_n)}\sup_{F_1\in\mathcal{F}_1,\dots,F_n\in\mathcal{F}_n}\rho_{h^*}(f(\mathbf{a},\mathbf{X})).
		\end{equation}
	\item	If the supremum of the right-hand side of \eqref{eq:eq-multi} is attained by some $F_1\in\mathcal{F}_1,\dots,F_n\in\mathcal{F}_n$ and $F\in\mathcal{D}(F_1,\dots,F_n)$, then for all $\mathbf{a}\in A$, $F_1^{\mathcal {I}_h},\dots,F_n^{\mathcal {I}_h}$ and a comonotonic random vector $({X}_1^{\I_h},\dots,{X}_n^{\I_h})$ with ${X}_i^{\I_h}\sim F_i^{\I_h}$, $i=1,\dots,n$ attain {the suprema} on both sides of \eqref{eq:eq-multi}.\footnote{A random vector $(X_1,\dots,X_n) \in (\L^1)^n$ is called \emph{comonotonic} if there exists a random variable $Z\in\X$ and increasing functions $f_1,\dots,f_n$ on $\R$ such that $X_i=f_i(Z)$ almost surely for all $i=1,\dots,n$.}
		\end{enumerate}
\end{theorem}	

Some examples of functions on $\R^n$ that are supermodular and positively homogeneous are given below. These functions are concave due to Theorem 3 of \cite{MM08}.
\begin{example}[Supermodular and positively homogeneous functions]
For $n\in\N$, the following functions $f:\R^n\to\R$ are supermodular and positively homogeneous. Write $\mathbf{x}=(x_1,\dots,x_n)\in\R^n$.
\begin{enumerate}[(i)]
\item (Linear function) $f:\mathbf{x}\mapsto\mathbf{a}^{\top}\mathbf{X}$ for $\mathbf{a}\in\R^n$. The function is increasing for $\mathbf{a}\in\R_+^n$.
\item (Geometric mean) $f:\mathbf{x}\mapsto-(\prod^n_{i=1}|x_i|)^{1/n}$ on $ \R_-^n$ for odd $n$. The function is also increasing on $\R_-^n$.
\item (Negated $p$-norm) $f:\mathbf{x}\mapsto -\Vert \mathbf{x}\Vert _p$  for $p\ge 1$. The function is increasing on $\R_-^n$.
\item (Sum of functions) $f:\mathbf{x}\mapsto\sum^n_{i=1}f_i(x_i)$ for positively homogeneous functions $f_1,\dots,f_n:\R\to\R$. The function is increasing if $f_1,\dots,f_n$ are increasing.
\end{enumerate}
\end{example}

	
	\section{One-dimensional uncertainty set with moment constraints}
	\label{sec:two_con}
	A popular example of an uncertainty set closed under concentration for all intervals is that of distributions with specified moment constraints as in Example \ref{ex:main}. We investigate this uncertainty set in detail and offer in this section some general results, which generalize several existing results in the literature; none of the results in the literature include non-monotone and non-convex distortion functions.  Non-monotone distortion functions create difficulties because of    possible complications at their discontinuity points.  
	
	For $p>1$, $m\in \R$ and $v>0$, we recall the set of interest in Example \ref{ex:main}:
	$$
	\M(p,m,v) = \{ F_Y\in \mathcal M_{p}: \E[Y]=m, ~\E[|Y-m|^{p}]\le v^{p}\}. 
	$$
	Let $q\in[1,\infty]$ be the H\"older conjugate of $p$, namely $q=(1-1/p)^{-1}$, or equivalently, $1/p+1/q=1$. For all $h\in \H^*$ or $h\in \H_*$, we denote by
	\begin{equation}
		\Vert h'-x\Vert _q= \left(\int_0^1 |h'(t)-x|^q \d t\right)^{1/q}, ~q<\infty
		\mbox{~~and~~}\Vert h'-x\Vert _\infty=\max_{t\in[0,1]}|h'(t)-x|,~~x\in\R. \label{eq:pnorm}
	\end{equation}
	We   introduce the following quantities:
	$$
	c_{h,q} =  \argmin_{x\in \R} \Vert  h'-x\Vert _q \mbox{~~~and~~~} [h]_q =  \min_{x\in \R} \Vert h'-x\Vert _q =  \Vert h'-c_{h,q}\Vert _q.$$
	We set  $[h]_q =\infty$ if $h$ is not continuous.
	It is easy to verify that $c_{h,q}$ is unique for $q>1$.
	The quantity $[h]_q$  may be interpreted as a $q$-central norm of the function $h$ and $c_{h,q}$ as its $q$-center. Note that 
	for $q=2$ and $h$ continuous, $[h]_2 = \Vert h'-h(1)\Vert _2$ and $c_{h,2}=h(1)$. We also note that the optimization problem is trivial if $[h]_q=0$, which corresponds to the case that $h' =h(1) \id_{[0,1]}$ and $\rho_h$ is a linear functional, thus a multiple of the expectation.
	In this case, the supremum and infimum are attained by all random variables whose distributions are in $\M(p,m,v)$, and they are equal to $m h(1)$. 
	Furthermore, for $h\in\mathcal H^*$ or $h\in \mathcal H_*$, and $q>1$, 
	we define a function on $[0,1]$ by
	$$
	\phi^q_h (t) =  \frac{|h'(1-t)-c_{h,q}| ^{q}}{h'(1-t)-c_{h,q}} [h]_q^{1-q}~~~ \mbox{if~}h'(1-t)-c_{h,q}\ne 0, \mbox{~~ and  $\phi^q_h (t)=0$ otherwise}. 
	$$ 
	In case $q=2$, for $t\in [0,1]$, $\phi^2_h(t)= (h'(1-t)-h(1))\Vert h'-h(1)\Vert _2^{-1}$ if $\Vert h'-h(1)\Vert _2>0$ and $0$ otherwise. We summarize our findings in the following theorem.
	\begin{theorem}\label{th:th1}
		For   any $h\in \mathcal H$, $m\in \R$, $v>0$ and $p>1$,  
		we have
		\begin{equation} \label{eq:th1} 
			\sup_{F_Y\in \mathcal M(p,m,v)} \rho_h(Y) =   m h(1) +  v [h^*]_q \mbox{~~and~~} \inf_{F_Y\in \mathcal M(p,m,v)} \rho_h(Y) =    m h(1)  -  v [h_*]_q.
		\end{equation}
		Moreover, if $h=\hat{h}$, $0<[h_*]_q<\infty$ and $0<[h^*]_q<\infty$, then  the  supremum and infimum in \eqref{eq:th1} are attained by a random variable $X$ such that $F_X\in \mathcal M(p,m,v)$ with  its  
		quantile function uniquely  specified as a.e.~equal to $ m+ v  \phi_{h^*}^q$ and 
		$ m-  v  \phi_{h_*}^q$, respectively. 
	\end{theorem} 
	
		The proof of Theorem \ref{th:th1} follows from a combination of Lemmas \ref{lem:p} and \ref{lem:p2} in Appendix \ref{appen:tech} and Theorem \ref{th:1}.
Note that 
for $h\in \mathcal H^*$ (resp.~$h\in \mathcal H_*$) and $q>1$, 
		$\phi^q_h $ is increasing (resp.~decreasing) on $[0,1]$.
		Hence, $\phi^q_h $ (resp.~$-\phi^q_h $) in Theorem \ref{th:th1}   indeed determines a quantile function. 

The following proposition concerns the finiteness of $\rho_h$ on $\mathcal L^p$.
	
	\begin{proposition}\label{co:finiteness}
		For any $h\in \H$ and $p\in[1,\infty]$, $\rho_h$ is finite on $\L^p$ if $[h^*]_q<\infty $ and $[h_*]_q<\infty$.
	\end{proposition} 
	
	As a special case of Proposition \ref{co:finiteness},  $\rho_h$ is always finite on $\L^1$ if $h$ is convex or concave  with bounded $h'$ because $[h^*]_\infty<\infty $ and $[h_*]_\infty<\infty$.  
	
	As a common example of the general result in Theorem \ref{th:th1}, below we collect our findings for the case of VaR.
	\begin{corollary}\label{prop:single-var}
		For $\alpha \in (0,1)$, $p>1$, $m\in \R$ and $v>0$, we have
		$$
		\sup_{F_Y\in \M(p,m,v)}\VaR_\alpha(Y)=\max_{F_Y\in \M(p,m,v)}\ES_\alpha(Y)= m+ v\alpha \left( \alpha^p (1-\alpha)+ (1-\alpha)^p \alpha\right)^{-1/p},
		$$
		and
		$$
		\inf_{F_Y\in \M(p,m,v)}\VaR_\alpha(Y)=\min_{F_Y\in \M(p,m,v)}\LE_\alpha(Y) =  m- v (1-\alpha) \left( \alpha^p (1-\alpha)+ (1-\alpha)^p \alpha\right)^{-1/p},
		$$
		where $$
		\LE_\alpha(Y)= \frac{1}{\alpha} \int_0^\alpha \VaR_t(Y)\d t, ~~Y\in \mathcal L^1.
		$$
	\end{corollary}


	We see from Theorem \ref{th:th1}  that if $h=\hat{h}$, then the supremum and the infimum of $\rho_{h}(Y)$ over $F_Y\in\M(p,m,v)$ are always attainable. 
	However, in case $h\neq \hat{h}$, the supremum or infimum may no longer be attainable as a maximum or minimum. We illustrate this in Example \ref{ex:1} below.


	\begin{example}[VaR and ES, $p=2$]
		\label{ex:1}
		Take $\alpha \in (0,1)$, $p=2$ and $\rho_h=\VaR_\alpha$, which implies $\rho_{h^*}=\ES_\alpha$.
		Corollary \ref{prop:single-var} gives
		$\sup_{F_Y\in \mathcal M(2,m,v)} \VaR_\alpha (Y) =  \sup_{F_Y\in \mathcal M(2,m,v)} \ES_\alpha (Y)  =   m +  v \sqrt{\alpha/(1-\alpha)}$.
		This is the well-known Cantalli-type formula for ES. 
		By Lemma \ref{lem:p}, the unique left-quantile function of the random variable $Z$ that attains the supremum of $\ES_\alpha$ is given by
		$F^{-1}_Z(t)= m + v (\id_{(\alpha,1]} (t)/(1-\alpha) -1) \sqrt{(1-\alpha)/\alpha}$, $t\in[0,1]$ a.e.
		We thus have $\VaR_\alpha(Z)=m - v \sqrt{ (1-\alpha)/(\alpha)}$, and hence $Z$ does not attain $ \sup_{F_Y\in \mathcal M(2,m,v)} \VaR_\alpha (Y)$. 
		It follows by the uniqueness of $F_Z$ that the supremum of 
		$ \VaR_\alpha (Y)$ over $F_Y\in \mathcal M(2,m,v)$  cannot be attained. However, the supremum of $\VaR^+_\alpha$ is attained by $Z$ since $\VaR^+_\alpha(Z)=m + v \sqrt{ (1-\alpha)/(\alpha)}$.
	\end{example}
	
	\begin{example}[Difference of two TK distortion riskmetrics]
	\label{ex:diff_con}
		Take $p=2$ and $h=h_1-h_2$ to be the difference between two inverse-S-shaped functions  in \eqref{eq:TKdistortion} with parameters the same as those in Example \ref{ex:diff} ($\gamma_1=0.8$, $\gamma_2=0.7$). 
		By Theorem \ref{th:th1}, the worst-case distortion riskmetrics under the uncertainty set $\M(2,m,v)$ are given by
		$\sup_{F_Y\in \mathcal M(2,m,v)}\rho_{h}(Y)=\sup_{F_Y\in \mathcal M(2,m,v)}\rho_{h^*}(Y)=0.3345v$,
		and the unique left-quantile function of the random variable $Z$ attaining both suprema above is given by
		$F^{-1}_Z(t)=m+2.9892\cdot h^{*\prime}(1-t)v$, $t\in[0,1]$ a.e.  
		The worst-case distortion riskmetrics obtained above are independent of the mean $m$ as $h(1)=h_1(1)-h_2(1)=0$, which is sensible since $\rho_h$ and $\rho_{h^*}$ only incorporate the disagreement between two distortion riskmetrics. Similarly, we can calculate the infimum of $\rho_h(Y)$ over $Y\in\M(2,m,v)$, and thus obtain the largest absolute difference between the two preferences numerically represented by $\rho_{h_1}$ and $\rho_{h_2}$. 
	\end{example}

	\section{Related optimization problems}
	\label{sec:app}
	
	In this section, we discuss the applications of our main results to some related optimization problems commonly investigated in the literature by including the outer problem of \eqref{eq:opt}.
		
	\subsection{Portfolio optimization}
	\label{sec:port}
	
	Our equivalence results can be applied to   robust portfolio optimization problems. For an uncertainty set $\widetilde{\M}\subset\M^n_p$ with $p\in[1,\infty]$, let the random vector $\mathbf{X}=(X_1,\dots,X_n)\sim F_{\mathbf{X}}\in\widetilde{\M}$, representing the random losses from $n$ risky assets. For $A\subset\R^n$, denote by a vector $\mathbf{a}\in A$ the  amounts invested in each of the $n$ risky assets. For a distortion function $h\in\H$ and distortion riskmetric $\rho_h:\L^p\to\R$, we aim to solve the robust portfolio optimization problem given by
	\begin{equation}
		\label{eq:port}
		\min_{\mathbf{a}\in A} \left( \sup_{F_{\mathbf{X}}\in\widetilde{\M}}\rho_h(\mathbf{a}^{\top}\mathbf{X})+\beta(\mathbf{a}) \right),
	\end{equation}
	where $\beta:\R^n\to\R$ is a penalty function of risk concentration. Note that  $\beta$ is irrelevant for the inner problem of \eqref{eq:port}. For a general non-concave $h$, there is no known algorithm to solve the inner problem of \eqref{eq:port}. 
	The outer optimization problem is also nontrivial in general. Therefore, we usually cannot obtain closed-form solutions of \eqref{eq:port} using classical results of optimization problems for non-convex risk measures. However, as a direct consequence of Theorems \ref{th:1} and \ref{th:multi}, the following proposition converts \eqref{eq:port} to an equivalent convex optimization problem that becomes much easier to solve. The proof of Proposition \ref{prop:port} follows directly from Theorems \ref{th:1} and \ref{th:multi}.
	\begin{proposition}
		\label{prop:port}
		Let  $h\in\H$, $n\in\N$, $A\subset\R^n$, and $\widetilde{\M}\subset\M^n_1$.
	\begin{enumerate}[(i)]
		\item if $h=\hat{h}$ and the set
		$
			\{F_{\mathbf{a}^{\top}\mathbf{X}}\in\M_1:F_{\mathbf{X}}\in\widetilde{\M}\}
		$
		is closed under concentration within $\I_h$ for all $\mathbf{a}\in A$, then 
		\begin{equation}
			\label{eq:port_equiv}
			\begin{aligned}
				\min_{\mathbf{a}\in A}\left(\sup_{F_{\mathbf{X}}\in\widetilde{\M}}\rho_h(\mathbf{a}^{\top}\mathbf{X})+\beta(\mathbf{a})\right)=\min_{\mathbf{a}\in A}\left(\sup_{F_{\mathbf{X}}\in\widetilde{\M}}\rho_{h^*}(\mathbf{a}^{\top}\mathbf{X})+\beta(\mathbf{a})\right).
			\end{aligned}
		\end{equation}
		\item if the set
		$
			\{F_{\mathbf{a}^{\top}\mathbf{X}}\in\M_1:F_{\mathbf{X}}\in\widetilde{\M}\}
		$
		is closed under concentration for all intervals for all $\mathbf{a}\in A$, then \eqref{eq:port_equiv} holds.
		
		\item If   $\widetilde{\mathcal M}$  is closed under 
		conditional expectation, then \eqref{eq:port_equiv} holds.
	\end{enumerate}
	\end{proposition}

	\subsection{Preference robust optimization}

We are also able to solve the preference robust optimization problem with distributional uncertainty. For $n\in\N$, an $n$-dimensional action set $A$, a set of plausible distributions $\widetilde{\M}\subset\M^n_1$, and a set of possible probability perceptions $\mathcal G\subset \mathcal H$, the problem is formulated as follows:
	\begin{equation}\label{eq:pro1}
		\min_{\mathbf{a}\in A} ~\sup_{F_\mathbf{X}\in \widetilde{\mathcal M} } ~\sup_{h\in \mathcal G}  \rho_h (f(\mathbf{a},\mathbf{X})).
	\end{equation}  
	Preference robust optimization refers to the situation when the objective is not completely known, e.g.,  $h$ is in the set $\mathcal G$ but not identified.
	Therefore, optimization is performed under the worst-case preference in the set $\mathcal G$. Also note that the form $\sup_{h\in \mathcal G}  \rho_{h} $ includes (but is not limited to) all coherent risk measures via the   representation of \cite{K01}. 
	For the problem of \eqref{eq:pro1} without distributional uncertainty (thus, only the minimum and the second supremum), see \cite{DL18}. 
	We have the following result whose proof follows from Theorems \ref{th:1} and \ref{th:multi}.
	\begin{proposition}
		\label{prop:pre}
		Let $\widetilde{\M}\subset\M^n_1$ and $A\subset\R^n$ with $n\in\N$.
	\begin{enumerate}[(i)]
		\item If $h=\hat{h}$ and the set 
		$
			\{F_{f(\mathbf{a},\mathbf{X})}\in\M_1:F_\mathbf{X}\in\widetilde{\M}\}
		$
		is closed under 
		concentration within $\I_h$ for all $\mathbf{a}\in A$,
		then for all $\mathcal G \subset \mathcal H$,
		\begin{equation}
			\label{eq:pre_equiv}
			\min_{\mathbf{a}\in A} \sup_{F_\mathbf{X}\in \widetilde{\M} } \sup_{h\in \mathcal G}  \rho_h (f(\mathbf{a},\mathbf{X})) = \min_{\mathbf{a}\in A} \sup_{F_\mathbf{X}\in \widetilde{\M} } \sup_{h\in \mathcal G}  \rho_{h^*} (f(\mathbf{a},\mathbf{X})).
		\end{equation}
		
		\item If the set 
		$
			\{F_{f(\mathbf{a},\mathbf{X})}\in\M_1:F_\mathbf{X}\in\widetilde{\M}\}
		$
		is closed under 
		concentration for all intervals for all $\mathbf{a}\in A$,
		then \eqref{eq:pre_equiv} holds for all $\mathcal G \subset \mathcal H$.
		
		\item If $\mathcal{G}$ is a set of increasing functions in $\H$, $f:A \times \R^n\to\R$ is concave in the second component,  and $\widetilde{\mathcal M}$  is closed under 
		conditional expectation, then \eqref{eq:pre_equiv} holds.
	\end{enumerate}
	\end{proposition}
	The preference robust optimization problem without distributional uncertainty (i.e., problem \eqref{eq:pro1} with only the minimum and the second supremum) is generally difficult to solve when the distortion function $h$ is not concave. However, when the distribution of the random variable is not completely known, we can transfer the original non-convex problem to its convex counterpart using \eqref{eq:pre_equiv}, provided that the set of plausible distributions is well structured.

	\section{Applications and numerical illustrations}
	\label{sec:app_num}
	
	Following the discussion in Section \ref{sec:app}, we provide several applications of our theoretical results to portfolio management   for specific sets of plausible distributions. None of the considered optimization problems in this section are convex, and we provide numerical calculations or approximation for the solutions to these optimization problems.\footnote{The processors we use are Intel(R) Xeon(R) CPU E5-2690 v3 @ 2.60GHz 2.59GHz (2 processors). The numerical results are calculated by MATLAB.}

	\subsection{Difference of risk measures under moment constraints}
	\label{sec:61}
	We demonstrate a price competition problem as an application of optimizing the difference between two risk measures shown in Example \ref{ex:diff_con}. Similar to the portfolio management problem discussed in Section \ref{sec:port}, we consider $n$ risky assets with random losses $X_1,\dots,X_n\in\L^2$ that are only known to have a fixed mean and a constrained covariance. That is, we choose the set $$\widetilde{\M}=\{F_{\mathbf{X}}\in\M^n_2:\E[\mathbf{X}]=\bm{\mu},~\mathrm{var}(\mathbf{X})\preceq\Sigma\},$$
	for $\bm{\mu}\in\R^n$ and $\Sigma\in \R^{n\times n}$ positive semidefinite. For an $n$-dimensional $\mathbf{a}\in A$, the set of all possible distributions of aggregate portfolio losses 
	\begin{equation}
	\label{eq:pre_meancov}
		\{F_{\mathbf{a}^{\top}\mathbf{X}}\in\M_2:F_{\mathbf{X}}\in\widetilde{\M}\}=\M^{\mathrm{mv}}(\mathbf{a},\bm{\mu},\Sigma)=\M\left(2,\mathbf{a}^{\top}\bm{\mu},\left(\mathbf{a}^{\top}\Sigma\mathbf{a}\right)^{1/2}\right)
	\end{equation}
	is closed under concentration for all intervals as is shown in Example \ref{ex:main}. Let $\rho_{h_1}:\L^2\to\R$ be an investor's own price of the portfolio, while $\rho_{h_2}:\L^2\to\R$ is her opponent's price of the same portfolio. We choose $h_1$ and $h_2$ to be the inverse-S-shaped distortion functions in \eqref{eq:TKdistortion}, with parameters the same as those in Example \ref{ex:diff_con} ($\gamma_1=0.8$ and $\gamma_2=0.7$). Write $h=h_1-h_2$. For an action set $A=\{(a_1,\dots,a_n)\in[0,1]^n:\sum^n_{i=1}a_i=1\}$,  the investor chooses the optimal $\mathbf{a}^*\in A$, such that the worst-case overpricing from her opponent is minimized.
	
	From the calculation in Example \ref{ex:diff_con}, we get
	\begin{equation}
	\label{eq:diff_opt}
	\begin{aligned}
		D(\Sigma)&:=\min_{\mathbf{a}\in A}\sup_{F_\mathbf{X}\in\widetilde{\M}}\left(\rho_{h_1}(\mathbf{a}^{\top}\mathbf{X})-\rho_{h_2}(\mathbf{a}^{\top}\mathbf{X})\right)\\
		&=\min_{\mathbf{a}\in A}\sup_{F_Y\in\M^{\mathrm{mv}}(\mathbf{a},\bm{\mu},\Sigma)}\rho_{h^*}(Y)=0.3345 \times \min_{\mathbf{a}\in A}\left(\mathbf{a}^{\top}\Sigma\mathbf{a}\right)^{1/2}.
	\end{aligned}
	\end{equation} 
	We note that optimizing $\rho_{h_1}-\rho_{h_2}$ is generally nontrivial since the difference between two distortion functions $h_1-h_2$ is not necessarily monotone, concave, or continuous, even though $h_1$ and $h_2$ themselves may have these properties. The generality of our equivalence result allows us to convert the original problem to the much simpler form \eqref{eq:diff_opt}, which can be solved efficiently.\footnote{The convex problem \eqref{eq:diff_opt} is solved numerically by the constrained nonlinear multivariable function ``fmincon" with the interior-point method.} Table \ref{tab:num_diff} demonstrates the optimal values of $\mathbf{a}^*$ and $D$ for different choices of $\Sigma$.
	\begin{table}[htbp]
		\centering
		\footnotesize
		\caption{Optimal results in \eqref{eq:diff_opt} for difference between two TK distortion riskmetrics}
		\begin{center}
		\begin{tabular}{cccc}
			$n$ & $\Sigma$ & $\mathbf{a}^*$ & $D$ \\\midrule
			$3$ & $\left(\begin{matrix}
			1 & 0 & 0\\0 & 1 & 0\\0 & 0 & 1
			\end{matrix}\right)$ & $(0.333,~0.333,~0.333)$ & $0.193$\\\midrule
			$3$ & $\left(\begin{matrix}
			2 & -1 & 0\\-1 & 2 & -1\\0 & -1 & 2
			\end{matrix}\right)$ & $(0.300,~0.400,~0.300)$ & $0.150$\\\midrule
			$3$ & $\left(\begin{matrix}
			1 & 1 & 1\\1 & 2 & 1\\1 & 1 & 3
			\end{matrix}\right)$ & $(0.997,~0.002,~0.001)$ & $0.335$\\\midrule
			$5$ & $\left(\begin{matrix}
			1 & 0 & 0 & 0 & 0\\0 & 2 & 0 & 0 & 0\\0 & 0 & 3 & 0 & 0\\0 & 0 & 0 & 4 & 0\\0 & 0 & 0 & 0 & 5
			\end{matrix}\right)$ & $(0.438,~0.219,~0.146,~0.110,~0.088)$ & $0.221$\\\bottomrule
		\end{tabular}
		\end{center}
		\label{tab:num_diff}
	\end{table}
	
	\subsection{Preference robust portfolio optimization with moment constraints}
	\label{sec:62}
	
	Next, we discuss an example of preference robust optimization with distributional uncertainty using the results in Sections \ref{sec:two_con}. Similarly to Section \ref{sec:61}, we consider the set of plausible aggregate portfolio loss distributions  
	\begin{equation*}
		\M^{\mathrm{mv}}(\mathbf{a},\bm{\mu},\Sigma)=\{F_{\mathbf{a}^{\top}\mathbf{X}}\in\M_2:F_{\mathbf{X}}\in\M^n_2,~\E[\mathbf{X}]=\bm{\mu},~\mathrm{var}(\mathbf{X})\preceq\Sigma\}
	\end{equation*}
	and the action set  $A=\{(a_1,\dots,a_n)\in[0,1]^n:\sum^n_{i=1}a_i=1\}$ representing the weights the investor assigns to each random loss.
	The investor considers TK distortion riskmetrics, however, she is not certain about the parameter $\gamma$ of the distortion function $h$. Thus, the investor consider the set of TK distortion riskmetrics  with distortion functions in
	$$\mathcal{G}=\left\{h\in\H:h=h^\gamma,~\gamma\in[0.5,0.9]\right\},$$ which is approximately the two-sigma confidence interval of $\gamma$ in \cite{WG96}.\footnote{The aggregate least-square estimate of $\gamma$ in Section 5 of \cite{WG96} is $0.71$ with standard deviation $0.1$.} 
    Therefore, the investor aims to find a optimal portfolio given the uncertainty in the riskmetrics. To penalize deviations from the benchmark parameter $\gamma=0.71$ \citep{WG96}, the investor use the term $ \mathrm{e}^{c(\gamma-0.71)^2}$ for some $c\ge 0$. 
	Since the set $\M^{\mathrm{mv}}(\mathbf{a},\bm{\mu},\Sigma)$ is closed under concentration for all intervals for all $\mathbf{a}\in  A$, Proposition \ref{prop:pre}, \eqref{eq:pre_meancov}, and Theorem \ref{th:th1} lead to 
	\begin{equation}
		\label{eq:pre_sol}
		\begin{aligned}
			V(\bm{\mu},\Sigma)&:=\min_{\mathbf{a}\in A}\sup_{F_Y\in\M^{\mathrm{mv}}(\mathbf{a},\bm{\mu},\Sigma)}\sup_{\gamma\in[0.5,0.9]}\left(\rho_{h^{\gamma}}(Y)-\mathrm{e}^{c(\gamma-0.71)^2}\right)\\&=\min_{\mathbf{a}\in A}\sup_{F_Y\in\M\left(2,\mathbf{a}^{\top}\bm{\mu},\left(\mathbf{a}^{\top}\Sigma\mathbf{a}\right)^{1/2}\right)}\sup_{\gamma\in[0.5,0.9]}\left(\rho_{(h^{\gamma})^*}(Y)-\mathrm{e}^{c(\gamma-0.71)^2}\right)\\&=\min_{\mathbf{a}\in A}\sup_{\gamma\in[0.5,0.9]}\left(\mathbf{a}^{\top}\bm{\mu}+\left(\mathbf{a}^{\top}\Sigma\mathbf{a}\right)^{1/2}\,[(h^{\gamma})^*]_2-\mathrm{e}^{c(\gamma-0.71)^2}\right).
		\end{aligned}
	\end{equation}
	
	We calculate the optimal values $V$ for different choices of parameters ($n$, $c$, $\bm{\mu}$ and $\Sigma$) and report them in Table \ref{tab:num_pre}, where $\mathbf{a}^*$ and $\hat\gamma$ represent the optimal weights and the parameters of the inverse-S-shaped distortion function, respectively. 
	Note that the last optimization problem in \eqref{eq:pre_sol} can be calculated numerically.\footnote{The problem \eqref{eq:pre_sol} is solved numerically by the constrained nonlinear multivariable function ``fmincon" with the interior-point method.}
	\begin{table}[htbp]
		\centering
		\footnotesize
		\caption{Optimal values in \eqref{eq:pre_sol} for TK distortion riskmetrics}
		\begin{center}
		\begin{tabular}{ccccccc}
			$n$ & $c$ & $\bm{\mu}$ & $\Sigma$ & $\mathbf{a}^*$ & $\hat\gamma$ & $V$ \\\midrule
			$3$ & $0$ & $(1,1,1)$ & $\left(\begin{matrix}
			1 & 0 & 0\\0 & 1 & 0\\0 & 0 & 1
			\end{matrix}\right)$ & $(0.333,~0.333,~0.333)$ & $0.610$ & $1.41$\\\midrule
			$3$ & $30$ & $(2,1,1)$ & $\left(\begin{matrix}
			1 & 0 & 0\\0 & 1 & 0\\0 & 0 & 1
			\end{matrix}\right)$ & $(0.000,~0.500,~0.500)$ & $0.676$ & $1.29$\\\midrule
			$3$ & $30$ & $(1,1,1)$ & $\left(\begin{matrix}
			2 & -1 & 0\\-1 & 2 & -1\\0 & -1 & 2
			\end{matrix}\right)$ & $(0.300,~0.400,~0.300)$ & $0.690$ & $1.17$\\\midrule
			$3$ & $30$ & $(1.2,1,1)$ & $\left(\begin{matrix}
			1 & 1 & 1\\1 & 2 & 1\\1 & 1 & 3
			\end{matrix}\right)$ & $(0.500,~0.331,~0.168)$ & $0.630$ & $1.57$\\\midrule
			$5$ & $30$ & $(1,1,1,1,1)$ & 
			$\left(\begin{matrix}
			1 & 0 & 0 & 0 & 0\\0 & 2 & 0 & 0 & 0\\0 & 0 & 3 & 0 & 0\\0 & 0 & 0 & 4 & 0\\0 & 0 & 0 & 0 & 5
			\end{matrix}\right)$
			 & $(0.438,~0.219,~0.146,~0.110,~0.088)$ & $0.678$ & $1.26$\\\bottomrule
		\end{tabular}
		\end{center}
		\label{tab:num_pre}
	\end{table}

	\subsection{Portfolio optimization with marginal constraints}
	\label{sec:63}
	A special case of the portfolio optimization problem introduced in Section \ref{sec:port}, which is of interest in robust risk aggregation (see e.g., \cite{BLLW20}), is to take $\widetilde{\M}$ to be the Fr\'echet class defined as \begin{equation}
		\label{eq:margin}
		\widetilde{\M}(F_1,\dots,F_n)=\{F_{\mathbf{X}}\in\M^n_1:X_i\sim F_i, ~i=1,\dots,n\},
	\end{equation}
	for some known marginal distributions $F_1,\dots,F_n\in\M_1$. In this case, although the left-hand side of \eqref{eq:port_equiv} is generally difficult to solve, for $A\subset\R^n_+$, the right-hand side of \eqref{eq:port_equiv} can be rewritten using convexity and comonotonicity as 
	\begin{equation}
	\label{eq:port_linear}
	\min_{\mathbf{a}\in A}\left(\mathbf{a}^{\top}({\rho}_{h^*}(X_1),\dots,{\rho}_{h^*}(X_n))+\beta(\mathbf{a})\right),
	\end{equation}
	where $X_i\sim F_i$, $i=1,\dots,n$. We see that \eqref{eq:port_linear} is a linear optimization problem with a penalty $\beta$, which often admits closed-form solutions when $\beta$ is properly chosen. For any given $\mathbf{a}\in A$, we define
	\begin{equation}
		\label{eq:weighted}
		\M(\mathbf{a},F_1,\dots,F_n)=\{F_{\mathbf{a}^{\top}\mathbf{X}}\in\M_1:X_i\sim F_i,~i=1,\dots,n\}.
	\end{equation}
	The set $\M(\mathbf{a},F_1,\dots,F_n)$ is the weighted version of $\M^{S}(F_1,\dots,F_n)$ in Example \ref{ex:main}. Note that $\M(\mathbf{a},F_1,\dots,F_n)$ is generally neither closed under concentration for all intervals nor closed under conditional expectation. However, $\M(\mathbf{a},F_1,\dots,F_n)$ is asymptotically (for large $n$) similar to a set of distributions closed under concentration for all intervals; see Theorem 3.5 of \cite{MW15} for a precise statement in the case of equal weights and identical marginal distributions. Therefore, even though $\M(\mathbf{a},F_1,\dots,F_n)$ is not closed under concentration for all intervals for some $\mathbf{a}\in A$, our result of the problem \eqref{eq:port_linear} is a good approximation of the original problem for large $n$. Such asymptotic equivalence between worst-case riskmetrics of aggregate risks with equal weights has already been well studied in the literature; see e.g., Theorem 3.3 of \cite{EWW15} for the $\VaR$/$\ES$ pair and Theorem 3.5 of \cite{CLW18} for distortion risk measures.

	We conduct numerical calculations to illustrate the equivalence between both sides in \eqref{eq:port_equiv}. 
    We choose the action set 
	$A_{a,b}=\{(x_1,\dots,x_n)\in[a,b]^n:\sum^n_{i=1}x_i=1\}$, for $0\le a<1/n<b\le 1$ and the penalty function $\beta$ to be the  $\mathcal{L}^2$-norm multiplied by a scaler $c\ge 0$, namely $c\Vert \cdot\Vert _2$, where the scaler $c$ is a tuning parameter of the $\mathcal L^2$ penalty. We first solve the optimization problems separately for the well-known VaR/ES pair at the level of $0.95$. Specifically, the two problems are given by
	\begin{align}
		V_{\mathrm{VaR}}(a,b,F_1,\dots,F_n)&=\min_{\mathbf{a}\in A_{a,b}}\left(\sup_{F_{\mathbf{X}}\in\M(F_1,\dots,F_n)}\VaR_{0.95}(\mathbf{a}^{\top}\mathbf{X})+c\Vert \mathbf{a}\Vert _2\right),\label{eq:var}\\
		V_{\mathrm{ES}}(a,b,F_1,\dots,F_n)&=\min_{\mathbf{a}\in A_{a,b}}\left(\sup_{F_{\mathbf{X}}\in\M(F_1,\dots,F_n)}\ES_{0.95}(\mathbf{a}^{\top}\mathbf{X})+c\Vert \mathbf{a}\Vert _2\right)\nonumber\\
		&=\min_{\mathbf{a}\in A_{a,b}}\left(\mathbf{a}^{\top}(\ES_{0.95}(F_1),\dots,\ES_{0.95}(F_n))+c\Vert \mathbf{a}\Vert _2\right),\label{eq:es}
	\end{align}
	where the true value of the original inner VaR problem is approximated by the rearrangement algorithm (RA) of \cite{PR12} and \cite{EPR13}, whereas the optimal value of the inner ES problem is obtained by simultaneously minimizing the sum of a linear combination of ES and the $2$-norm of the vector $\mathbf{a}$, which can be done efficiently.\footnote{The outer problems of \eqref{eq:var} and \eqref{eq:es} are solved numerically by the constrained nonlinear multivariable function ``fmincon" with the sequential quadratic programming (SQP) algorithm. The same method is also applied when solving outer problems of \eqref{eq:non-convex} and \eqref{eq:convex}.} In particular, if the marginals of the random losses are identical (i.e., $F_1=\cdots=F_n=F$), the optimal solution is $\mathbf{a^*}=(1/n,\dots,1/n)$ and $V_{\mathrm{ES}}(a,b,F_1,\dots,F_n) = \ES_{0.95}(F)+c/\sqrt{n}$. We consider the following marginal distributions
	\begin{enumerate}[(i)]
	\item $F_i$ follows a Pareto distribution with scale parameter $1$ and shape parameter $3+(i-1)/(n-1)$ for $i=1,\dots,n$;
	\item $F_i$ is normally distributed with parameters $\mathrm{N}(1,1+(i-1)/(n-1))$, for $i=1,\dots,n$;  
	\item $F_i$ follows an exponential distribution with parameter $1+(i-1)/(n-1)$, for $i=1,\dots,n$. 
	\end{enumerate} We choose $n$ to be $3$, $10$, and $20$. For comparison, we calculate the value $n\Vert \Delta\mathbf{a}^*\Vert _2$, where $\Delta\mathbf{a}^*$ is the difference between the optimal weights of the non-convex problem and the convex problem. In addition, we calculate the absolute differences between the optimal values obtained by the two problems, $\Delta V  =V_{\mathrm{ES}}-V_{\mathrm{VaR}}\ge 0$,  and the percentage differences $\Delta V /V_{\mathrm{VaR}}$. Tables \ref{tab:num_var1} and \ref{tab:num_var2} show the numerical results that compare both optimization problems with two choices of the action sets $A_{a,b}$. The computation time is reported (in seconds).  We observe that the optimal values obtained in the two problems get closer and become approximately the same as $n$ gets larger. As explained before, this is because the set of plausible distributions $\M(F_1,\dots,F_n)$ is asymptotically equal to a set  closed under concentration for all intervals.

	Next, we consider a TK distortion riskmetric with parameter $\gamma=0.7$. 
	Due to the non-concavity of $h$, there are no known   ways of directly solving the non-convex optimization problem 
	\begin{equation}
		\label{eq:non-convex}
		\min_{\mathbf{a}\in A_{a,b}}\left(\sup_{F_{\mathbf{X}}\in\M(F_1,\dots,F_n)}\rho_h(\mathbf{a}^{\top}\mathbf{X})+c\Vert \mathbf{a}\Vert _2\right).
	\end{equation} 
	We may get an approximation of \eqref{eq:non-convex} using a lower bound of  $\rho_h$ in \eqref{eq:non-convex} produced with the dependence structure created by the rearrangement algorithm (RA);\footnote{Such a dependence structure  obviously provides a lower bound for the worst-case value in \eqref{eq:non-convex}. In theory, the result from RA is thus not an optimal dependence structure for \eqref{eq:non-convex}. In our numerical results,   this lower bound is very close to an upper bound only for the case of VaR and ES but not for the case of TK distortion riskmetrics.} for simplicity, we denote   by $V_h$ this lower bound.
	 On the other hand, by \eqref{eq:port_equiv}, the convex counterpart of \eqref{eq:non-convex} can be written (using Theorem \ref{th:1}) as
	\begin{equation}
		\label{eq:convex}
		\begin{aligned}
			V_{h^*}(a,b,F_1,\dots,F_n)&=\min_{\mathbf{a}\in A_{a,b}}\left(\sup_{F_{\mathbf{X}}\in\M(F_1,\dots,F_n)}\rho_{h^*}(\mathbf{a}^{\top}\mathbf{X})+c\Vert \mathbf{a}\Vert _2\right)\\			&=\min_{\mathbf{a}\in A_{a,b}}\left(\mathbf{a}^{\top}({\rho}_{h^*}(X_1),\dots,{\rho}_{h^*}(X_n))+c\Vert \mathbf{a}\Vert _2\right),
		\end{aligned}
	\end{equation}
where $X_i\sim F_i$ for $i=1,\dots,n$. We calculate the absolute differences between the optimal values of the convex and non-convex problems $\Delta V=V_{h^*}-V_h\ge 0$, and the percentage differences $\Delta V/V_h$. Tables \ref{tab:numerical1} and \ref{tab:numerical2} compare the numerical results of the two optimization problems with different choices of $A_{a,b}$. We observe that the percentage differences between the RA lower bound $V_h$ for the non-convex problem \eqref{eq:non-convex}  and the minimum value $V_{h^*}$  of the convex problem \eqref{eq:convex} are roughly between $10\%$ to $20\%$. Note that  the RA lower bound is not expected to be very close to the true minimum of \eqref{eq:non-convex}, and hence the differences between the solution of \eqref{eq:non-convex} and the optimal value in \eqref{eq:convex} are smaller than the observed numbers. 


Note that, by transforming a non-convex optimization problem to a convex one, we significantly reduce the computational time of calculating bounds with negligible errors, as shown in Tables \ref{tab:num_var1}-\ref{tab:numerical2}.

	 \begin{table}[ht]
		\centering
		\small
		\caption{Comparison of the numerical results of the two optimization problems \eqref{eq:var}  and \eqref{eq:es}  for $\VaR_{0.95}$ and $\ES_{0.95}$ with $a=0$ and $b=1$}
		\begin{center}
		\begin{tabular}{c c c|cc|cc| c c c}
						 & &  $ c$ & $V_{\mathrm{VaR}}$ & time & $V_{\mathrm{ES}}$  & time &    \multirow{1}{*}{$n\Vert \Delta\mathbf{a}^*\Vert _2$} & \multirow{1}{*}{$\Delta V $} & \multirow{1}{*}{$ {\Delta V }/{V_{\mathrm{VaR}}}$ (\%)}\\  
			\midrule
			\multirow{3}{*}{\begin{minipage}{0.06\textwidth} \begin{center} (i)\\ Pareto \end{center}\end{minipage}} & $n=\phantom{0}3$ & $2.5$ & $3.547$ & $\phantom{0}31.53$ & $3.741$ & $0.72$ & $8.88\times 10^{-5}$ & $0.194\phantom{0}$ & $5.48\phantom{0}$\\
			& $n=10$ & $3.0$ & $3.197$ & $153.83$ & $3.215$ & $1.39$ & $9.18\times 10^{-4}$ & $0.0178$ & $0.558$\\
			& $n=20$ & $4.0$ & $3.156$ & $424.17$ & $3.159$ & $9.37$ & $3.53\times 10^{-5}$ & $2.68\times 10^{-3}$ & $\phantom{0}0.0850$ \\
			\midrule
			\multirow{3}{*}{\begin{minipage}{0.06\textwidth} \begin{center} (ii)\\ Normal \end{center}\end{minipage}} & $n=\phantom{0}3$ & $4.0$ & $5.766$ & $\phantom{0}31.19$ & $5.785$ & $0.18$ & $1.39\times 10^{-3}$ & $0.0186$ & $0.323$\\
			& $n=10$ & $2.0$ & $4.082$ & $\phantom{0}97.30$ & $4.083$ & $0.77$ & $1.18\times 10^{-3}$ & $3.24\times 10^{-5}$ & $7.93\times 10^{-4}$\\
			& $n=20$ & $3.0$ & $4.132$ & $431.79$ & $4.132$ & $4.66$ & $2.69\times 10^{-3}$ & $1.88\times 10^{-5}$ & $4.55\times 10^{-4}$\\\midrule
			\multirow{3}{*}{\begin{minipage}{0.06\textwidth} \begin{center} (iii)\\ Exp \end{center}\end{minipage}} & $n=\phantom{0}3$ & $3.0$ & $4.251$ & $\phantom{0}26.78$ & $4.405$ & $0.07$ & $0.331$ & $0.155\phantom{0}$ & $3.64\phantom{0}$\\
			& $n=10$ & $4.0$ & $3.892$ & $118.23$ & $3.893$ & $0.50$ & $9.74\times 10^{-4}$ & $2.92\times 10^{-4}$ & $7.52\times 10^{-3}$\\
			& $n=20$ & $7.0$ & $4.230$ & $543.03$ & $4.230$ & $3.47$ & $3.08\times 10^{-4}$ & $4.47\times 10^{-5}$ & $1.06\times 10^{-3}$\\\bottomrule
		\end{tabular}
		\end{center}
		\label{tab:num_var1}
	\end{table} 
	
	\begin{table}[htbp]
		\centering
		\small
		\caption{Comparison of the numerical results of the two optimization problems \eqref{eq:var}  and \eqref{eq:es} for $\VaR_{0.95}$ and $\ES_{0.95}$ with $a=1/(2n)$ and $b=2/n$}
		\begin{center}
		\begin{tabular}{ccc|cc|cc|ccc}
						 & &  $ c$ & $V_{\mathrm{VaR}}$ & time & $V_{\mathrm{ES}}$  & time &    \multirow{1}{*}{$n\Vert \Delta\mathbf{a}^*\Vert _2$} & \multirow{1}{*}{$\Delta V $} & \multirow{1}{*}{$ {\Delta V }/{V_{\mathrm{VaR}}}$ (\%)}\\  
\midrule
			\multirow{3}{*}{\begin{minipage}{0.06\textwidth} \begin{center} (i)\\ Pareto \end{center}\end{minipage}} & $n=\phantom{0}3$ & $2.5$ & $3.546$ & $\phantom{0}54.59$ & $3.741$ & $0.19$ & $6.58\times 10^{-4}$ & $0.194\phantom{0}$ & $5.48\phantom{0}$\\
			& $n=10$ & $3.0$ & $3.204$ & $146.63$ & $3.220$ & $1.60$ & $1.99\times 10^{-4}$ & $0.0160$ & $0.498$\\
			& $n=20$ & $4.0$ & $3.162$ & $847.13$ & $3.163$ & $10.08\phantom{0}$ & $1.69\times 10^{-3}$ & $2.23\times 10^{-3}$ & $\phantom{0}0.0706$\\\midrule
			\multirow{3}{*}{\begin{minipage}{0.06\textwidth} \begin{center} (ii)\\ Normal \end{center}\end{minipage}} & $n=\phantom{0}3$ & $4.0$ & $5.766$ & $\phantom{0}57.31$ & $5.785$ & $0.19$ & $1.32\times 10^{-3}$ & $0.0187$ & $0.324$\\
			& $n=10$ & $2.0$ & $4.084$ & $166.25$ & $4.084$ & $0.79$ & $0$ & $2.94\times 10^{-5}$ & $7.20\times 10^{-4}$\\
			& $n=20$ & $3.0$ & $4.133$ & $691.91$ & $4.133$ & $5.91$ & $0$ & $1.99\times 10^{-5}$ & $4.82\times 10^{-4}$\\\midrule
			\multirow{3}{*}{\begin{minipage}{0.06\textwidth} \begin{center} (iii)\\ Exp \end{center}\end{minipage}} & $n=\phantom{0}3$ & $3.0$ & $4.369$ & $\phantom{0}48.58$ & $4.422$ & $0.09$ & $1.04\times 10^{-3}$ & $0.0533$ & $1.22\phantom{0}$\\
			& $n=10$ & $4.0$ & $3.916$ & $115.18$ & $3.916$ & $0.50$ & $2.54\times 10^{-5}$ & $1.38\times 10^{-4}$ & $3.52\times 10^{-3}$\\
			& $n=20$ & $7.0$ & $4.236$ & $665.05$ & $4.236$ & $3.48$ & $2.73\times 10^{-4}$ & $4.04\times 10^{-5}$ & $9.54\times 10^{-4}$\\\bottomrule
		\end{tabular}
		\end{center}
		\label{tab:num_var2}
	\end{table}

	\begin{table}[ht]
		\centering
		\small
		\caption{Comparison of the numerical results of the two optimization problems \eqref{eq:non-convex} and \eqref{eq:convex} for TK distortion riskmetrics with $a=0$ and $b=1$}
		\begin{center}
		\begin{tabular}{ccc|cc|cc|ccc} 
	& &  $ c$ & $V_h$ & time & $V_{h^*}$  & time &    \multirow{1}{*}{$n\Vert \Delta\mathbf{a}^*\Vert _2$} & 	\multirow{1}{*}{$\Delta V $} & \multirow{1}{*}{$ {\Delta V }/{V_{h}}$ (\%)}\\  
			\midrule
			\multirow{3}{*}{\begin{minipage}{0.06\textwidth} \begin{center} (i)\\ Pareto \end{center}\end{minipage}} & $n=\phantom{0}3$ & $1.0$ & $1.076$ & $144.75$ & $1.185$ & $0.23$ & $0.488$ & $0.109$ & $10.2$\\
			& $n=10$ & $2.0$ & $1.047$ & $220.03$ & $1.237$ & $1.42$ & $0$ & $0.190$ & $18.1$\\
			& $n=20$ & $4.0$ & $1.301$ & $826.64$ & $1.501$ & $8.24$ & $0$ & $0.200$ & $15.4$\\\midrule
			\multirow{3}{*}{\begin{minipage}{0.06\textwidth} \begin{center} (ii)\\ Normal \end{center}\end{minipage}} & $n=\phantom{0}3$ & $0.5$ & $1.240$ & $\phantom{0}60.76$ & $1.493$ & $0.16$ & $\phantom{0}0.0784$ & $0.253$ & $20.4$\\
			& $n=10$ & $0.5$ & $1.141$ & $246.31$ & $1.363$ & $0.72$ & $1.28\phantom{0}$ & $0.222$ & $19.4$\\
			& $n=20$ & $0.5$ & $1.103$ & $1503.35\phantom{0}$ & $1.316$ & $2.80$ & $1.78\phantom{0}$ & $0.213$ & $19.3$\\\midrule
			\multirow{3}{*}{\begin{minipage}{0.06\textwidth} \begin{center} (iii)\\ Exp \end{center}\end{minipage}} & $n=\phantom{0}3$ & $1.0$ & $1.305$ & $\phantom{0}49.79$ & $1.427$ & $0.23$ & $0.360$ & $0.122$ & $9.32$\\
			& $n=10$ & $2.0$ & $1.313$ & $198.43$ & $1.484$ & $1.62$ & $0.184$ & $0.171$ & $13.0$\\
			& $n=20$ & $2.0$ & $1.120$ & $850.12$ & $1.286$ & $10.91\phantom{0}$ & $0.158$ & $0.166$ & $14.8$\\\bottomrule 
		\end{tabular}
		\end{center}
		\label{tab:numerical1}
	\end{table}

	\begin{table}[ht]
		\centering
		\small
		\caption{Comparison of the numerical results of the two optimization problems \eqref{eq:non-convex} and \eqref{eq:convex} for TK distortion riskmetrics with $a=1/(2n)$ and $b=2/n$}
		\begin{center}
		\begin{tabular}{ccc|cc|cc|ccc} 
		& &  $ c$ & $V_h$ & time & $V_{h^*}$  & time &    \multirow{1}{*}{$n\Vert \Delta\mathbf{a}^*\Vert _2$} & 	\multirow{1}{*}{$\Delta V $} & \multirow{1}{*}{$ {\Delta V }/{V_{h}}$ (\%)}\\
			\midrule
			\multirow{3}{*}{\begin{minipage}{0.06\textwidth} \begin{center} (i)\\ Pareto \end{center}\end{minipage}} & $n=\phantom{0}3$ & $1.0$& $1.077$ & $\phantom{0}73.21$ & $1.185$ & $0.25$ & $0.469$ & $0.109$ & $\phantom{0}10.11$\\
			& $n=10$ & $2.0$& $1.047$ & $248.38$ & $1.237$ & $2.29$ & $0.378$ & $0.191$ & $18.2$\\
			& $n=20$ & $4.0$ & $1.301$ & $638.24$ & $1.501$ & $12.21\phantom{0}$ & $0$ & $0.200$ & $15.4$\\\midrule
			\multirow{3}{*}{\begin{minipage}{0.06\textwidth} \begin{center} (ii)\\ Normal \end{center}\end{minipage}} & $n=\phantom{0}3$ & $0.5$ & $1.240$ & $179.68$ & $1.493$ & $0.19$ & $\phantom{0}0.0784$ & $0.253$ & $20.4$\\
			& $n=10$ & $0.5$ & $1.146$ & $389.97$ & $1.363$ & $0.76$ & $0.660$ & $0.217$ & $19.0$\\
			& $n=20$ & $0.5$ & $1.103$ & $1563.84\phantom{0}$ & $1.316$ & $3.39$ & $1.63\phantom{0}$ & $0.213$ & $19.3$\\\midrule
			\multirow{3}{*}{\begin{minipage}{0.06\textwidth} \begin{center} (iii)\\ Exp \end{center}\end{minipage}} & $n=\phantom{0}3$ & $1.0$ & $1.304$ & $\phantom{0}52.66$ & $1.430$ & $0.25$ & $0.107$ & $0.126$ & $9.65$\\
			& $n=10$ & $2.0$ & $1.312$ & $236.15$ & $1.485$ & $2.27$ & $0.214$ & $0.172$ & $13.1$\\
			& $n=20$ & $2.0$ & $1.119$ & $879.73$ & $1.289$ & $10.10\phantom{0}$ & $0.141$ & $0.170$ & $15.2$\\\bottomrule 
		\end{tabular}
		\end{center}
		\label{tab:numerical2}
	\end{table}

	\section{Concluding remarks} 
	\label{sec:conclude}
	We introduced the new concept of closedness under concentration, which is, in the context of distributional uncertainty, a sufficient condition  to transform an optimization problem with a non-convex distortion riskmetric  to its convex counterpart. 
	This concept is genuinely weaker than closedness under conditional expectation,
	and our main result unifies and improves many existing results in the literature. 
	Many sets of plausible distributions commonly used in the literature of finance, optimization, and risk management are closed under concentration within some $\I$.  Moreover, by focusing on distortion riskmetrics whose distortion functions are not necessarily monotone, concave, or continuous, we are able to solve optimization problems for the class of functionals larger than classical risk measures or deviation measures. In particular, we are able to obtain  bounds on  differences between two distortion riskmetrics, which represent measures of disagreement between two utilities/risk attitudes.
	Our result can also be applied to solve the popular problem of optimizing  risk measures under moment constraints. In particular, we obtain the worst- and best-case distortion riskmetrics when the underlying random variable has a fixed mean and bounded $p$-th moment.  

	We demonstrate the applicability of our result by numerically calculating the solution to optimizing the difference between risk measures, preference robust optimization  and  portfolio optimization under marginal constraints. In all numerical examples, the original non-convex problem is converted or well approximated by a convex one which can be solved efficiently.  
	
%
	
Our condition of closedness under concentration within $\I$ in Theorem \ref{th:1}  is  sufficient but not necessary   for the equivalence of a non-convex and a convex optimization problem under distributional uncertainty. A necessary  condition of the equivalence is closedness under concentration of the set of maximizers in Theorem \ref{th:nec}.
An open question is to find  
a \emph{necessary and sufficient} condition on the uncertainty set $\M$ itself such that the desired  equivalence holds. 
Pinning down such a condition may facilitate many more applications in  decision theory, finance, game theory, and operations research.

\subsection*{Acknowledgments}
  The authors would like to thank anonymous referees for their constructive comments enhancing the paper. SMP would like to acknowledge the support of the Natural Sciences and Engineering Research Council of Canada with funding reference numbers DGECR-2020-00333 and RGPIN-2020-04289. RW acknowledges financial support from the Natural Sciences and Engineering Research Council of Canada (RGPIN-2018-03823, RGPAS-2018-522590).

\newpage
	\appendix
	

\normalsize
     \setcounter{lemma}{0}
     \renewcommand{\thelemma}{A.\arabic{lemma}} 
          \setcounter{proposition}{0}
     \renewcommand{\theproposition}{A.\arabic{proposition}}
          \setcounter{definition}{0}
     \renewcommand{\thedefinition}{A.\arabic{definition}}
               \setcounter{corollary}{0}
     \renewcommand{\thecorollary}{A.\arabic{corollary}}
               \setcounter{example}{0}
     \renewcommand{\theexample}{A.\arabic{example}} 
               \setcounter{equation}{0}
     \renewcommand{\theequation}{A.\arabic{equation}} 
     
               \setcounter{figure}{0}
     \renewcommand{\thefigure}{A.\arabic{figure}}

     \begin{center}
    \LARGE Technical appendices
     \end{center}
     
\section{Omitted technical details from the paper}
In this appendix, we present    technical details for some examples and as well as some technical remarks omitted from the paper.
\subsection{Proofs of   claims in some Examples}

\label{app:main}
%
\begin{proof}[Proof of the claim in Example \ref{ex:main}]
We show that $\M^{\mathrm{mv}}(\mathbf{a},\bm{\mu},\Sigma)$ is equivalent to
			$$\{F_S\in\M_2:\E[S]=\mathbf{a}^{\top}\bm{\mu},~\mathrm{var}(S)\le\mathbf{a}^{\top}\Sigma\mathbf{a}\}=\M\left(2,\mathbf{a}^{\top}\bm{\mu},\left(\mathbf{a}^{\top}\Sigma\mathbf{a}\right)^{1/2}\right).$$
			For a proof of the equivalence between the sets with fixed mean and covariance matrix, see \cite{P07}. Indeed, it is clear that $\M^{\mathrm{mv}}(\mathbf{a},\bm{\mu},\Sigma)\subset\M(2,\mathbf{a}^{\top}\bm{\mu},(\mathbf{a}^{\top}\Sigma\mathbf{a})^{1/2})$. On the other hand, for all $F_S\in\M(2,\mathbf{a}^{\top}\bm{\mu},(\mathbf{a}^{\top}\Sigma\mathbf{a})^{1/2})$, we write $\mathbf{a}=(a_1,\dots,a_n)$, $\bm{\mu}=(\mu_1,\dots,\mu_n)$, and take $\mathbf{X}=(X_1,\dots,X_n)$ such that $X_i=(S-\mathbf{a}^{\top}\bm{\mu})/(na_i)+\mu_i$, for $i=1,\dots,n$. It follows that $F_S=F_{\mathbf{a}^{\top}\mathbf{X}}\in\M^{\mathrm{mv}}(\mathbf{a},\bm{\mu},\Sigma)$. Therefore, we have $\M^{\mathrm{mv}}(\mathbf{a},\bm{\mu},\Sigma)=\M(2,\mathbf{a}^{\top}\bm{\mu},(\mathbf{a}^{\top}\Sigma\mathbf{a})^{1/2})$.
\end{proof}
\begin{proof}[Proof of the claim in Example \ref{ex:wassball}]
We will show that $\M(\widetilde{G},\epsilon)$ is closed under concentration within $\I$ for all $\I\subset\widetilde{\I}$. Write $\I=\{C_i:i\in K\}$ for some $K\subset\N$. For all $i\in K$ and $F\in\M(\widetilde{G},\epsilon)$, we have $\widetilde{G}^{-1}(u)=c_i$ for $u\in C_i$ for some $c_i\in\R$. For all $i\in K$, by Jensen's inequality, $$\frac{1}{\lambda(C_i)}\int_{C_i}\left|F^{-1}(u)-\widetilde{G}^{-1}(u)\right|^p\d u\ge \left|\frac{\int_{C_i}F^{-1}(u)\d u}{\lambda(C_i)}-c_i\right|^p=\frac{1}{\lambda(C_i)}\int_{C_i}\left|(F^{C_i})^{-1}(u)-\widetilde{G}^{-1}(u)\right|^p\d u.$$
It follows that
$$\begin{aligned}
(W_p(F,\widetilde{G}))^p-(W_p(F^{C_i},\widetilde{G}))^p&=\int^1_0\left|F^{-1}(u)-\widetilde{G}^{-1}(u)\right|^p\d u-\int^1_0\left|(F^{C_i})^{-1}(u)-\widetilde{G}^{-1}(u)\right|^p\d u\\
&=\int_{C_i}\left|F^{-1}(u)-\widetilde{G}^{-1}(u)\right|^p\d u-\int_{C_i}\left|(F^{C_i})^{-1}(u)-\widetilde{G}^{-1}(u)\right|^p\d u\ge 0,
\end{aligned}$$
and thus $W_p(F^{C_i},\widetilde{G})\le W_p(F,\widetilde{G})\le \epsilon$.
Moreover,   \eqref{eq:invG3} and the above argument lead to 
$$
(W_p(F,\widetilde{G}))^p-(W_p(F^{\I},\widetilde{G}))^p = \sum_{i\in  K} (W_p(F,\widetilde{G}))^p-(W_p(F^{C_i},\widetilde{G}))^p \ge 0.
$$
Hence, $W_p(F^{\mathcal I},\widetilde{G})\le W_p(F,\widetilde{G})\le \epsilon$. 
\end{proof}
\begin{proof}[Proof of the claim in Example \ref{ex:11}]
\label{app:11}
For $\epsilon\ge 0$, $\mathbf{w}\in[0,\infty)^n$, $p>1$, $a>1$ and $\mathbf{Z}\in(\L^p)^n$, by Theorem 7 of \cite{MWW22}, the uncertainty set
$$\{F_{\mathbf{w}^{\top}\mathbf{X}}\in\M_p:F_\mathbf{X}\in\mathcal{M}^n(F_{\mathbf{Z}},\epsilon)\}=\M(F_{\mathbf{w}^{\top}\mathbf{Z}},\epsilon\Vert \mathbf{w}\Vert _{b}),$$
where $b$ is the conjugate of $a$ (i.e., $1/a+1/b=1$). Suppose that for a benchmark distribution $\widetilde{G}\in\M^n_p$, there exists a random vector $\mathbf{Z}\sim\widetilde{G}$ such that $\mathbf{Z}\ge\mathbf 0$ and $\P(\mathbf{Z}=\mathbf 0)=p_0$ for some $p_0\in(0,1]$.
Note that $\p(\mathbf w^\top \mathbf Z =0) \ge p_0$ and the quantile function of $\mathbf w^\top \mathbf Z$ is equal to $0$ on $(0,p_0]$. 
 It follows from Example \ref{ex:wassball} that the set $\M(F_{\mathbf{w}^{\top}\mathbf{Z}},\epsilon\Vert \mathbf{w}\Vert _{b})$ is closed under concentration within $\{(0,t)\}$ for all $t\le p_0$.  
\end{proof}
\begin{proof}[{Proof of the claim in Example \ref{ex:opt_hedge}}]

We will show that the set of distributions,
$$\mathcal{M}=\{F_{X-V(X)+g(\E[V(X)])}\in\M_1:V\in\mathcal{V}\},$$
is closed under concentration within $\{(p,1)\}$ for all $p\in[p_0,1)$. For each $V\in\mathcal{V}$ and a standard uniform random variable $U$, we write $a=\E[F^{-1}_{X-V(X)}(U)|U\in(p,1)]$. Since $F^{-1}_X(p)\ge l$, we can take $$W(x)=V(x)\id_{\{x\le F^{-1}_X(p)\}}+(x-a)\id_{\{x>F^{-1}_X(p)\}},~~x\in\R.$$ It follows that $W\in\mathcal{V}$. Noting that $a=\E[X-V(X)|X>F^{-1}_X(p)]$, we have
$$\begin{aligned}
&X-W(X)+g(\E[W(X)])\\&=(X-V(X))\id_{\{X\le F^{-1}_X(p)\}}+a\id_{\{X>F^{-1}_X(p)\}}+g\left(\E[V(X)\id_{\{X\le F^{-1}_X(p)\}}+(X-a)\id_{\{X>F^{-1}_X(p)\}}]\right)\\
&=(X-V(X))\id_{\{X\le F^{-1}_X(p)\}}+a\id_{\{X>F^{-1}_X(p)\}}+g(\E[V(X)]),
\end{aligned}$$
which follows the same distribution as $F_{X-V(X)+g(\E[V(X)])}^{(p,1)}$.
It follows that $\mathcal{M}$ is closed under concentration within $\{(p,1)\}$ for all $p\in[p_0,1)$.
\end{proof}

\subsection{A few additional technical remarks mentioned in the paper}
\label{app:a5}
\begin{remark}[on Theorem \ref{th:1}]
\label{rem:infinite}		Using Theorem \ref{th:1}, if  for some $\mathbf{a}\in A$, the set $\M:=\{F_{f(\mathbf{a},\mathbf{X})}:F_\mathbf{X}\in\widetilde{\M}\}$ is closed under concentration for all intervals and $ \sup\{   \rho_{h^*}(f(\mathbf{a},\mathbf{X})): F_\mathbf{X}\in \widetilde{\mathcal M}\} =\infty$,
		then 
		$\sup\{   \rho_{h}(f(\mathbf{a},\mathbf{X})): F_\mathbf{X}\in \widetilde{\mathcal M}\} =\infty$. Thus, both objectives in the inner optimization of \eqref{eq:opt} are infinite for this $\mathbf a$, which can be excluded from the outer optimization over $ A$.  Verifying  $\sup\{   \rho_{h^*}(f(\mathbf{a},\mathbf{X})): F_\mathbf{X}\in \widetilde{\mathcal M}\} =\infty$ is easier than  verifying $ \sup\{   \rho_{h}(f(\mathbf{a},\mathbf{X})): F_\mathbf{X}\in \widetilde{\mathcal M}\} =\infty$ since generally $\rho_h$ is smaller than $\rho_{h^*}$.
	\end{remark}
	\begin{remark}[on Example \ref{ex:main}]\label{rem:cxorder}
	Using Strassen's Theorem (e.g., Theorem 3.A.4 of \cite{SS07}), closedness under conditional expectation can equivalently be expressed using convex order. A set $\mathcal M\subset \mathcal M_1$   is closed under conditional expectation if and only if it holds that for $F\in \mathcal M$ and $G\le_{\rm cx} F$, we have $G\in \mathcal M$.
	\end{remark}
	
	\begin{remark}[on  Proposition \ref{lem:lem-4}]
		\label{rem:counterex}
		In Proposition \ref{lem:lem-4}, if $\M$ is closed under conditional expectation, $\mathcal{I}$ can be taken as an infinite set. However, $\M$ may not be closed under concentration within an infinite $\I$ if  we only assume that $\M$ is closed under concentration for all intervals. Indeed, if we take $\M$ as the set of distributions obtained by some $F\in\M$ with finitely many concentrations, then clearly $\M$ is closed under concentration for all intervals. However, $F^{\mathcal{I}}\notin\M$ when $\mathcal{I}$ is an infinite collection of disjoint intervals. This also serves as a counter-example of the converse statement of Proposition \ref{lem:lem-2} since $\M$ is closed under concentration for all intervals but not closed under conditional expectation.
	\end{remark}
	
		\section{Proofs of all technical results}
	\label{app:proofs}
	
	We present all proofs of technical results in this appendix.	Throughout, we denote the set of discontinuity points of $h$ (excluding $0$ and $1$) by
	\begin{equation}\label{eq:discontinuity}
		J_h=\{t\in(0,1):h(t)\ne h(t^+) \mbox{ or } h(t)\ne h(t^-)\}.
	\end{equation}
	Note that $\hat h(t)$ can be written as 
			\begin{equation}
			\label{eq:hat_h}
			\hat{h}(t)=\left\{\begin{array}{l l}
				h(t^+)\vee h(t^-)\vee h(t), & t\in J_h,\\[0.25em]
				h(t), & \text{otherwise}.
			\end{array}\right.
		\end{equation} 
	
	\subsection{Proof  of results in Section \ref{sec:2}}

%
	
	\begin{proof}[Proof of Proposition \ref{lem:extension}]
		Note that $(\hat{h})^*=h^*=\hat{h}=h$ on $0$ and $1$. For all $t\in(0,1)$, since $(\hat{h})^*(t)\ge \hat{h}(t)\ge h(t)$, we have $(\hat{h})^*(t)\ge h^*(t)$. On the other hand, we have $h^*(t)\ge h(t^+)$ for $t\in(0,1)$. Indeed, if $h^*(t_0)<h(t^+_0)$ for some $t_0\in(0,1)$, then we have $h^*(t_0+\epsilon)<h(t_0+\epsilon)$ for some $\epsilon>0$, which leads to a contradiction. Similarly, we have $h^*(t)\ge h(t^-)$ for $t\in(0,1)$. Together with $h^*\ge h$ on $(0,1)$, we have $h^*\ge \hat{h}$ on $(0,1)$, which implies that $h^*\ge(\hat{h})^*$ on $(0,1)$. Therefore, $(\hat{h})^*=h^*$ on $[0,1]$. 
		
		Next, we assert that the set $\{t\in [0,1]: \hat{h}(t)\ne h^*(t)\}$ is a union of disjoint sets that are not singletons. To show this assertion,  assume that the converse is true. There exists $x\in(0,1)$, such that $\hat{h}(x)<h^*(x)$ and $\hat{h}(t)=h^*(t)$ on $t\in(x-\epsilon,x)\cup(x,x+\epsilon)$ for some $0<\epsilon\le x\wedge (1-x)$. It is clear that $x\in J_h$. Since $h^*$ is continuous on $(x-\epsilon,x+\epsilon)$, we have$$\hat{h}(x)<h^*(x)=h^*(x^+)=\hat{h}(x^+).$$ This contradicts   \eqref{eq:hat_h}. Therefore, the set $\{t\in [0,1]: \hat{h}(t)\ne h^*(t)\}$ is the union of some disjoint intervals, denoted by $\cup_{l\in L}A_l$ for some $L\subset\N$. For all $l\in L$, we denote the left and right endpoints of $A_l$ by $a_l$ and $b_l$, respectively, with $a_l<b_l$. 
		Define a function via linear interpolation 
		$$h^c(t)=\left\{\begin{array}{l l}
		\hat{h}(a_l)+\frac{\hat{h}(b_l)-\hat{h}(a_l)}{b_l-a_l}(t-a_l), & t\in A_l,~l\in L,\\
		\hat{h}(t), & \text{otherwise}.
		\end{array}\right.$$
		It is clear that $h^c\le h^*$ and $h^c$ is continuous on $(0,1)$. We will prove that $h^c=h^*$ on $\cup_{l\in L}A_l$. Suppose for the purpose of contradiction that $h^c\ne h^*$ on $\cup_{l\in L}A_l$. 
		Since $h^c<h^*$ for some point in  $\cup_{l\in L}A_l$, there exists $x_0\in A_l$ for some $l\in L$ such that $h^c(x_0)<\hat{h}(x_0)$. Thus we can take a point $(x_1,\hat{h}(x_1))\in (0,1)\times \R$ with $\hat{h}(x_1)>h^c(x_1)$, which has the largest perpendicular distance to the straight line $h^c(t)=\hat{h}(a_l)+\frac{\hat{h}(b_l)-\hat{h}(a_l)}{b_l-a_l}(t-a_l)$, namely,
		$$x_1=\argmax_{\substack{x\in A_l\\\hat{h}(x)>h^c(x)}}\frac{(b_l-a_l)\hat{h}(x)-(\hat{h}(b_l)-\hat{h}(a_l))x-(b_l-a_l)\hat{h}(a_l)+(\hat{h}(b_l)-\hat{h}(a_l))a_l}{ \left((\hat{h}(b_l)-\hat{h}(a_l))^2+(b_l-a_l)^2\right)^{1/2}}.$$
		The existence of the maximizer $x_1$ is due to the upper semicontinuity of $\hat h$.
		There exists a function $g$ with $g=h^*$ on $[0,1]\setminus A_l$ and $g(x_1)=\hat{h}(x_1)$, such that $g$ is concave and $\hat{h}\le g\le h^*$ on $[0,1]$. Since $h^*>\hat{h}$ on $A_l$, we have $h^*(x_1)>\hat{h}(x_1)=g(x_1)$. Thus $h^*$ cannot be the concave envelope of $\hat{h}$, which leads to a contradiction. Thus, $h^*=h^c$ on $\cup_{l\in L}A_l$. 
		Since $h^*=\hat h=h^c$ on $(0,1) \setminus  (\cup_{l\in L}A_l)$, we have $h^*=h^c$. 
		Therefore, $\{t\in [0,1]: \hat{h}(t)\ne h^*(t)\}$	is a union of disjoint open intervals, and $h^*$ is linear on each of the intervals.
%
	\end{proof}

	\subsection{Proofs  of results in Section \ref{sec:main}} 
	
	\begin{proof}[Proof of Theorem \ref{th:1}]
	We will first show that, assuming that $\M$  is closed under concentration within $\I_h$, we have 
		\begin{equation}\label{eq:result1}
			\sup_{F_X\in\M} \rho_{\hat{h}}(X)=\sup_{F_X\in\M} \rho_{h^*}(X).
		\end{equation}
	After proving \eqref{eq:result1}, we show the three statements in Theorem \ref{th:1} in the order {(i)}, {(ii)}, and {(iii)}. 
		
		 For $h\in\H$, suppose that $\M$ is closed under concentration within $\I_h$. Take an arbitrary random variable $Y$ with  
		$F_Y\in \mathcal M$. Let $G=F^{\mathcal{I}_h}_Y$. For $h\in \H$, write  functions $g(t)=1-\hat{h}(1-t)$ and $g_*(t)=1-h^*(1-t)$ for $t\in [0,1]$. By definition of $\mathcal{I}_h$, $g\ne g_*$ on each set in $\mathcal{I}_h$ and $g=g_*$ on other sets. For any $(a,b)\in \mathcal{I}_h$, we have
		$G^{-1} (t)= \frac{ \int_{a}^b F^{-1}_Y(u) \d u }{b-a}$ for all $t\in(a,b]$ and $G^{-1+} (t)= \frac{ \int_{a}^b F^{-1}_Y(u) \d u }{b-a}$ for all $t\in[a,b)$.
		Using the fact  that $g_*$ is linear on $ (a,b)$ and $g(t)=g_*(t)$ for $t=a,b$,  we have
		\begin{equation}
			\label{eq:invG2}
			\begin{aligned}
				\int_{(a,b)} F^{-1}_Y(t) \d {g_*}(t) &= ( g_*(b)-  g_*(a) )
				\frac{ \int_{a}^b F^{-1}_Y(t) \d t }{b-a} 
				\\&=( g(b)-  g(a) )
				\frac{ \int_{a}^b F^{-1}_Y(t) \d t }{b-a}\\
				&= \int_{(a,b]}G^{-1} ( t)  \d g(t)+G^{-1+}(a)(g(a^+)-g(a)).
			\end{aligned}
		\end{equation}
		Define the sets
		\begin{gather*}
			J_+=\{t\in J_h:\hat{h}(t^+)=\hat{h}(t)\neq\hat{h}(t^-)\},~~J_-=\{t\in J_h:\hat{h}(t^+)\neq\hat{h}(t)= \hat{h}(t^-)\},\\
			\text{and}~~J_0=\{t\in J_h:\hat{h}(t^+)\neq \hat{h}(t)\neq \hat{h}(t^-)\}.
		\end{gather*} 
		To better understand these sets, we recall Figure \ref{fig:hat-h} (without concave envelopes) as Figure \ref{fig:A1} to demonstrate an example of a distortion function $h$, the corrresponding $\hat{h}$, the sets $J_h$, $J_+$, $J_-$, and $J_0$, and the sets $\hat{J}$, $\hat{J}_+$, $\hat{J}_-$, $\hat{J}^0_+$, and $\hat{J}^0_-$ (defined in the proof of {(i)} below).
		 
		\begin{figure}[htbp]
		\begin{center}
			\includegraphics[width=0.45\textwidth]{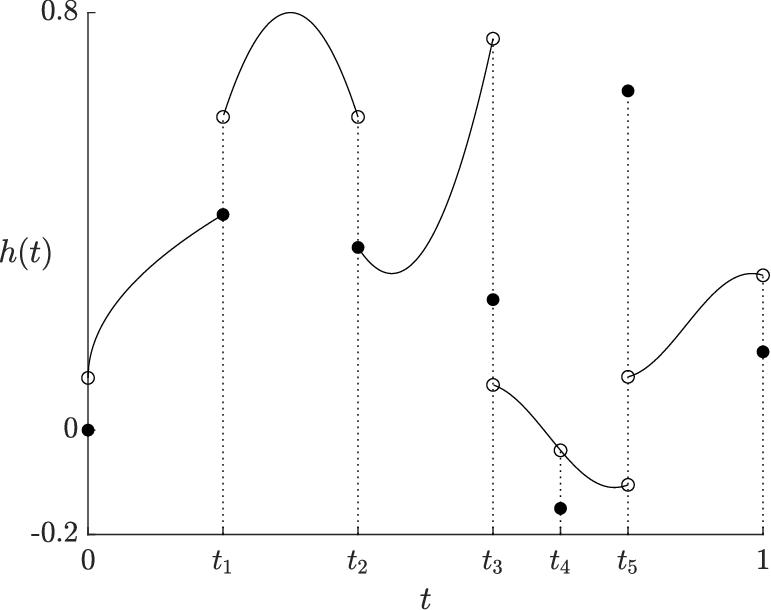}~~~~
			\includegraphics[width=0.45\textwidth]{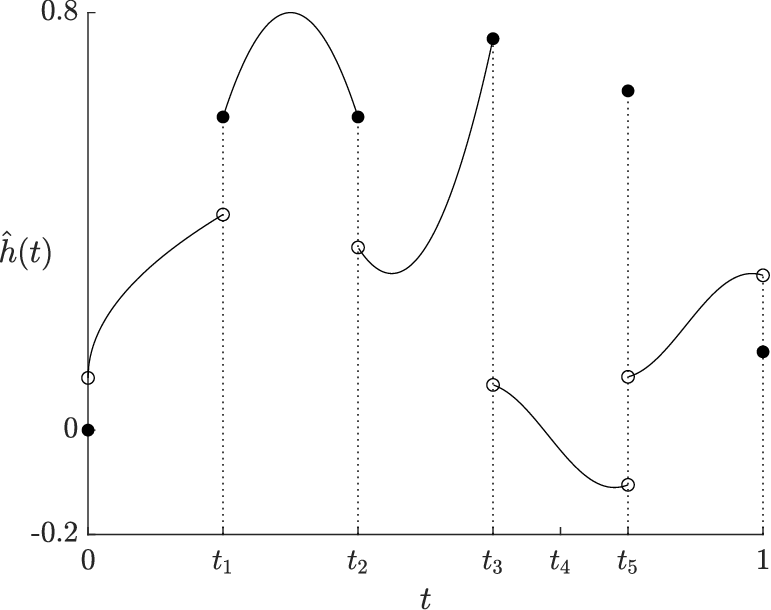}
		\end{center}
		\caption{An example of $h$ (left) and $\hat{h}$ (right); in this figure, $J_h=\{t_1,t_2,t_3,t_4,t_5\}$, $J_+=\{t_1\}$, $J_-=\{t_2,t_3\}$, and $J_0=\{t_5\}$. Moreover, the sets we use in the proof of (i) are $\hat{J}=\{t_1,t_2,t_3,t_4\}$, $\hat{J}_+=\{t_1,t_4\}$, $\hat{J}_-=\{t_2,t_3\}$, $\hat{J}^0_+=\{t_4\}$, and $\hat{J}^0_-=\{t_3\}$}
		\label{fig:A1}
		\end{figure}
		
		Note that for a random variable $Z_{\I_h}\sim F^{\I_h}_Y$, we have
		\begin{align*}
			\rho_{\hat{h}}(Z_{\mathcal{I}_h})&=\int_{(0,1]\setminus (J_+\cup J_0)} G^{-1}(t)\d g(t)+\sum_{t\in J_+\cup J_0\cup\{0\}}G^{-1+}(t)(g(t^+)-g(t)).
		\end{align*}
		Hence using \eqref{eq:invG2} and \eqref{eq:invG3}, we get
		\begin{equation}
			\label{eq:equiv}
			\begin{aligned}
				&\rho_{h^*}(Y)  -\rho_{\hat{h}}(Z_{\mathcal{I}_h}) \\ &= \int_0^1 F^{-1}_Y(t) \d {g_*}(t) + F^{-1+}_Y(0)(g_*(0^+)-g_*(0))  \\&\quad \quad -   \int_{(0,1]\setminus (J_+\cup J_0)} G^{-1}(t)\d g(t)-\sum_{t\in J_+\cup J_0\cup\{0\}}G^{-1+}(t)(g(t^+)-g(t))\\ 
				&= \sum_{(a,b) \in \mathcal{I}_h} \left(\int_{(a,b)} F^{-1}_Y(t) \d {g_*}(t) -\int_{(a,b]}G^{-1} ( t)  \d g(t)-G^{-1+}(a)(g(a^+)-g(a)) \right)=0.
			\end{aligned}
		\end{equation}		
		Since $\mathcal{M}$ is closed under concentration within $\I_h$, we have $F^{\mathcal{I}_h}_Y\in\M$ by definition. Thus we have
		\begin{equation*}
			\rho_{h^*}(Y)=\rho_{\hat{h}}(Z_{\mathcal{I}_h})\le\sup_{F_X\in\M} \rho_{\hat{h}}(X),
		\end{equation*}
		which gives our desired equality \eqref{eq:result1} since $\rho_{h^*}=\rho_{(\hat{h})^*} \ge\rho_{\hat{h}}$.
		
		Proof of {(i)}: Using $h=\hat{h}$ and \eqref{eq:result1}, we have $\sup_{F_X\in\M} \rho_{h}(X)=\sup_{F_X\in\M} \rho_{h^*}(X)$.

		Proof of {(ii)}:
		We will prove (ii) in two main steps. First, we show that (ii) holds if $\I_h$ is finite and $h$ has finitely many discontinuity points. Next, we discuss general $h$.
		
		\textbf{\emph{Finite case:}}	Here we prove \eqref{eq:eq-main} under the case where $\I_h$ is finite and $h$ has finitely many discontinuity points (i.e. $J_h$ in \eqref{eq:discontinuity} is a finite set). Suppose that $\M$ is closed under concentration for all intervals, it directly implies that $\M$ is closed under concentration within $\I_h$ by Proposition \ref{lem:lem-4}. Therefore,
		\eqref{eq:result1} holds for all $h\in\H$.
		 Next, we need to show that 
		$ 	\sup_{F_X\in\M} \rho_{{h}}(X)=\sup_{F_X\in\M} \rho_{\hat{h}}(X).$ 
		Define
		\begin{gather*}
			\hat{J}=\{t\in J_h: \hat{h}(t)\ne h(t)\},~~
			\hat{J}_+=\{t\in \hat{J}:\hat{h}(t)=\hat{h}(t^+)\},~~\text{and}~~\hat{J}_-=\hat{J}\setminus\hat{J}_+.
		\end{gather*}
		For $n>0$, write intervals
		$$A^n_s=\left\{\begin{array}{l l}
		(1-s-1/\sqrt{n},1-s+1/n), & s\in\hat{J}_-,\\
		(1-s-1/n,1-s+1/\sqrt{n}), & s\in \hat{J}_+.
		\end{array}\right.$$
		Let $\mathcal{I}^n=\{A^n_s:s\in\hat{J}\}$. Note that $h\in\H$ has finitely many discontinuity points. Thus the intervals in $\mathcal{I}^n$ are disjoint when $n$ is large enough. For all $F_Y\in\M$ and $Y\sim F_Y$, we define
		$$Z_{\mathcal{I}^n}=F^{-1}_Y(U)\id_{\{U\notin\bigcup_{s\in \hat{J}}A^n_s\}}+\sum_{s\in \hat{J}}\E[F^{-1}_Y(U)|U\in A^n_s]\id_{\{U\in A^n_s\}}.$$
		It follows that $Z_{\mathcal{I}^n}\sim F^{\mathcal{I}^n}_Y$ and the right-quantile function of $Z_{\mathcal{I}^n}$, denoted by $G^{-1+}_n$, is given by the right-continuous adjusted version of
		$$F^{-1+}_Y(t)\id_{\{t\notin\bigcup_{s\in \hat{J}}A^n_s\}}+\sum_{s\in \hat{J}}\frac{\int_{A^n_s}F^{-1}_Y(u)\d u}{\lambda(A^n_s)}\id_{\{t\in A^n_s\}},~~t\in(0,1).$$
		Thus we get
		$$\lim_{n\to\infty}G^{-1+}_n(1-t)=\left\{\begin{array}{l l}
		F^{-1}_Y(1-t), & t\in\hat{J}_-,\\
		F^{-1+}_Y(1-t), & \text{otherwise}.
		\end{array}\right.$$
		Similarly, if we denote the left-quantile function of $Z_{\mathcal{I}^n}$ by $G^{-1}_n$, then $G^{-1}_n$ is given by the left-continuous version of
		$$F^{-1}_Y(t)\id_{\{t\notin\bigcup_{s\in \hat{J}}A^n_s\}}+\sum_{s\in \hat{J}}\frac{\int_{A^n_s}F^{-1}_Y(u)\d u}{\lambda(A^n_s)}\id_{\{t\in A^n_s\}}.$$
		It follows that
		$$\lim_{n\to\infty}G^{-1}_n(1-t)=\left\{\begin{array}{l l}
		F^{-1+}_Y(1-t), & t\in\hat{J}_+,\\
		F^{-1}_Y(1-t), & \text{otherwise}.
		\end{array}\right.$$
		Define, further, the sets
		$$\hat{J}^0_+=\{t\in\hat{J}_+:h(t)\neq h(t^-)\}~~\text{and}~~\hat{J}^0_-=\{t\in\hat{J}_-:h(t)\neq h(t^+)\}.$$
		For $u\in[0,1]$, write
		$$\begin{aligned}
		&h_-(u)=\sum_{t\in\hat{J}_-}(h(t)-h(t^-))\id_{\{u\ge t\}},~~h^0_-(u)=\sum_{t\in\hat{J}^0_-}(h(t^+)-h(t))\id_{\{u> t\}},\\
		&h_+(u)=\sum_{t\in\hat{J}_+}(h(t^+)-h(t))\id_{\{u>t\}},~~h^0_+(u)=\sum_{t\in\hat{J}^0_+}(h(t)-h(t^-))\id_{\{u\ge t\}},\\
		&\hat{h}_-(u)=\sum_{t\in\hat{J}_-}(h(t^+)-h(t^-))\id_{\{u> t\}},~~\hat{h}_+(u)=\sum_{t\in\hat{J}_+}(h(t^+)-h(t^-))\id_{\{u\ge t\}},\\
		& \text{and}~~h_0(u)=h(u)-h_+(u)-h_-(u)-h^0_+(u)-h^0_-(u)=\hat{h}(u)-\hat{h}_+(u)-\hat{h}_-(u).
		\end{aligned}$$
		Note that $|Z_{\mathcal{I}^n}-F^{-1}_Y(U)|=0$ when $U\notin\bigcup_{s\in \hat{J}}A^n_s$ and $0,1\in [0,1]\setminus\bigcup_{s\in \hat{J}}A^n_s$. We have $|Z_{\mathcal{I}^n}-F^{-1}_Y(U)|<\infty$. Therefore, by the dominated convergence theorem,
		$$\begin{aligned}
		&\lim_{n\to\infty}(\rho_{h_-}(Z_{\mathcal{I}^n})+\rho_{h^0_-}(Z_{\mathcal{I}^n}))\\
		&=\lim_{n\to\infty}\int^1_0G^{-1+}_n(1-u)\d {h_-}(u)+\lim_{n\to\infty}\int^1_0G^{-1}_n(1-u)\d {h^0_-}(u)\\
		&=\sum_{t\in\hat{J}_-}F^{-1}_Y(1-t)(h(t)-h(t^-))+\sum_{t\in\hat{J}^0_-}F^{-1}_Y(1-t)(h(t^+)-h(t))\\
		&=\sum_{t\in\hat{J}_-\setminus\hat{J}^0_-}F^{-1}_Y(1-t)(h(t)-h(t^-))+\sum_{t\in\hat{J}^0_-}F^{-1}_Y(1-t)(h(t)-h(t^-)+h(t^+)-h(t))\\
		&=\sum_{t\in\hat{J}_-\setminus\hat{J}^0_-}F^{-1}_Y(1-t)(h(t^+)-h(t^-))+\sum_{t\in\hat{J}^0_-}F^{-1}_Y(1-t)(h(t^+)-h(t^-))=\rho_{\hat{h}_-}(Y).
		\end{aligned}$$
		Similarly, we get $\lim_{n\to\infty}(\rho_{h_+}(Z_{\mathcal{I}^n})+\rho_{h^0_+}(Z_{\mathcal{I}^n}))=\rho_{\hat{h}_+}(Y)$. On the other hand, it is clear that $\lim_{n\to\infty}\rho_{h_0}(Z_{\mathcal{I}^n})=\rho_{h_0}(Y)$. Therefore, we have
		$$\begin{aligned}
		\lim_{n\to\infty}\rho_{h}(Z_{\mathcal{I}^n})&=\lim_{n\to\infty}(\rho_{h_-}(Z_{\mathcal{I}^n})+\rho_{h^0_-}(Z_{\mathcal{I}^n})+\rho_{h_+}(Z_{\mathcal{I}^n})+\rho_{h^0_+}(Z_{\mathcal{I}^n})+\rho_{h_0}(Z_{\mathcal{I}^n}))\\
		&=\rho_{\hat{h}_-}(Y)+\rho_{\hat{h}_+}(Y)+\rho_{h_0}(Y)=\rho_{\hat{h}}(Y).
		\end{aligned}$$
		Thus we have
		\begin{equation}
			\label{eq:result3}
			\rho_{\hat{h}}(Y)=\lim_{n\to\infty}\rho_{h}(Z_{\mathcal{I}^n})\le\sup_{F_X\in\M}\rho_h(X).
		\end{equation}
		Using \eqref{eq:result1} and \eqref{eq:result3}, we get
		$$\sup_{F_X\in\M}\rho_{h^*}(X)=\sup_{F_X\in\M}\rho_{\hat{h}}(X)\le\sup_{F_X\in\M}\rho_h(X).$$

		\textbf{\emph{General case:}} We prove Theorem \ref{th:1} for all general $h\in\H$ where $\mathcal{I}_h$ or the number of discontinuity points of $h$ is countable.
		
			1. If $\mathcal{I}_h$ is countable, it suffices to prove \eqref{eq:result1}. We write $\mathcal{I}_h$ as the collection of $(a_i,b_i)$ for $i\in\N$ and define $\mathcal{I}^n_1=\{(a_i,b_i):i=1,\dots,n\}$ for all $n\in\N$. Define the function
			$$h_n(t)=\left\{\begin{array}{l l}
			h^*(t), & t\in(1-b_i,1-a_i),~i=1,\dots,n,\\
			\hat{h}(t), & \text{otherwise}.
			\end{array}\right.$$
			It is clear that for all $n\in\N$, the set $\{t\in[0,1]:h_n(t)\ne \hat{h}(t)\}$ is a finite union of disjoint open intervals and $h_n$ is linear on each of the intervals. For all random variables $Y$ with $F_Y\in\M$, let random variable $Z_{\I^n_1}\sim F^{\I^n_1}_Y$. Similar to \eqref{eq:result1}, we have
			$$\rho_{h_n}(Y)=\rho_{\hat{h}}(Z_{\mathcal {I}^n_1})\le\sup_{F_X\in\M} \rho_{\hat{h}}(X),~~\text{for all }n\in\N.$$
			Note that $h_n(t)\uparrow h^*(t)$ as $n\to\infty$ for all $t\in(0,1)$. By the monotone convergence theorem, we get $\rho_{h_n}(Y)\to\rho_{h^*}(Y)$ as $n\to\infty$. It follows that
			\begin{equation*}
				\sup_{F_X\in\M} \rho_{\hat{h}}(X)\ge\rho_{h_n}(Y)\xrightarrow{n\to\infty}\rho_{h^*}(Y).
			\end{equation*}
			
			2. If $h\in\H$ has countably many discontinuity points, it suffices to prove \eqref{eq:result3}. There exist series of finite sets $\{\hat{J}^m\}_{m\in\N}\subset\hat{J}$, such that $\hat{J}^m\to\hat{J}$ as $m\to \infty$. For all $m\in\N$, write
			$$\hat{h}_m(t)=\left\{\begin{array}{l l}
			\hat{h}(t), & t\in \hat{J}^m,\\
			h(t), & \text{otherwise},
			\end{array}\right.$$
			and define
			\begin{gather*}
				\hat{J}^m_+=\{t\in \hat{J}^m:\hat{h}_m(t)=\hat{h}_m(t^+)\},~~\text{and}~~\hat{J}^m_-=\hat{J}^m\setminus\hat{J}^m_+.
			\end{gather*}
			For $n>0$, let $\mathcal{I}^{n,m}_2=\{B^{n,m}_{s}:i\in\hat{J}^m\}$ with
			$$B^{n,m}_{s}=\left\{\begin{array}{l l}
			(1-s-1/\sqrt{n},1-s+1/n), & s\in\hat{J}^m_-,\\
			(1-s-1/n,1-s+1/\sqrt{n}), & s\in \hat{J}^m_+.
			\end{array}\right.$$
			Following the same argument as \eqref{eq:result3}, for all random variable $Y$ with $F_Y\in\M$, we have
			$$\sup_{F_X\in\M}\rho_h(X)\ge\rho_{h}(Z_{\mathcal{I}^{n,m}_2})\xrightarrow{n\to\infty}\rho_{\hat{h}_m}(Y),~~\text{for all }m\in\N,$$
			where $Z_{\I^{n,m}_2}\sim F^{\I^{n,m}_2}_Y$.
			Moreover, we have $\hat{h}_m(t)\uparrow\hat{h}(t)$ for all $t\in[0,1]$ as $m\to\infty$. By the monotone convergence theorem, we have $\rho_{\hat{h}_m}(Y)\to\rho_{\hat{h}}(Y)$ as $m\to\infty$. Therefore, we have
			\begin{equation*}
				\sup_{F_X\in\M}\rho_{\hat{h}}(X)\le\sup_{F_X\in\M}\rho_h(X).
			\end{equation*}

		Proof of {(iii)}: For all $h\in\H$, if $\mathcal M$ is closed under concentration within $\I_h$ and $h=\hat{h}$, we have $F_Y^{\mathcal {I}_h}\in \mathcal M$ by definition. Since $Z_{\mathcal{I}_h}\sim F_Y^{\mathcal {I}_h}$, \eqref{eq:equiv} gives
		$$\rho_{h^*}(Y)   = \rho_{\hat{h}}(Z_{\mathcal{I}_h})=\rho_h(Z_{\mathcal{I}_h}).$$ Note that $\rho_h\le \rho_{h^*}$ generally. 
		Therefore,  if 
		$ \max_{F_Y\in \mathcal M}  \rho_{h^*}(Y) 
		$ is attained by   $F_Y$, then so is 
		$ \max_{F_Y\in \mathcal M } \rho_h(Y) $ by $ F_Y^{\mathcal {I}_h}$.
		Obviously, these two quantities share a common maximizer $F_Y^{\mathcal {I}_h}$ because 
		$$ \rho_{h^*}(Z_{\mathcal{I}_h})\le  \max_{F_Y\in \mathcal M}  \rho_{h^*}(Y) =   \max_{F_Y\in \mathcal M}  \rho_{h }(Y)   =  \rho_h(Z_{\mathcal{I}_h})\le \rho_{h^*}(Z_{\mathcal{I}_h}).$$
		The proof is complete.
	\end{proof}

	\begin{proof}[Proof of Theorem \ref{th:nec}]
%
Suppose for contradiction that $\M_{\mathrm{opt}}$ is not closed under concentration within $\I_h$. There exists $F_Y\in\M_{\mathrm{opt}}$, such that $F^{\I_h}_Y\notin\M_{\mathrm{opt}}$. Define the set $\mathcal{Y}_h=\{(F^{-1}_Y(a),F^{-1}_Y(b)):(a,b)\in\I_h\}$. Since $F^{\I_h}_Y\notin\M_{\mathrm{opt}}$, there exists an interval $(a,b)\in\I_h$, such that $F^{-1}_Y$ is not constant on $(a,b)$. Thus the Lebesgue measure $\lambda((F^{-1}_Y(a),F^{-1}_Y(b)))>0$. Since $h^*>h$ on $(a,b)$,
\begin{equation}\label{eq:difference}
\begin{aligned}
\rho_{h^*}(Y)-\rho_{h}(Y)&=\int_\R(h^*(\P(Y>x))-h(\P(Y>x)))\d x\\&=\sum_{A\in\mathcal{Y}_h}\int_A(h^*(\P(Y>x))-h(\P(Y>x)))\d x>0.
\end{aligned}
\end{equation}
On the other hand, we have
\begin{equation*}
\rho_{h^*}(Y)\le \sup_{F_X\in \mathcal M}  \rho_{h^*}(X)=\sup_{F_X\in \mathcal M}  \rho_{h}(X)=\rho_h(Y)\le\rho_{h^*}(Y),
\end{equation*}
which leads to a contradiction to \eqref{eq:difference}. Therefore, $\mathcal{M}_{\mathrm{opt}}$ is closed under concentration within $\I_h$. 
\end{proof}

\begin{proof}[Proof of Proposition \ref{lem:lem-2}]
		We first prove that closedness under conditional expectation implies closedness under concentration for all intervals. For all random variables $Y\in\L^1$ and intervals $C\subset [0,1]$, define
		$$X=F^{-1}_Y(U)\id_{\{U\not \in C\}} +  \E[F^{-1}_Y(U)|U\in C]\id_{\{U  \in C\}},$$
		where $U\sim\mathrm{U}[0,1]$.
		The distribution of $X$ is the concentration $F^C_Y$. For all $\sigma(X)$-measurable random variables $Z$, we have that $Z|\{U\in C\}$ is constant. Hence,
		$$\begin{aligned}
		\E[XZ]&=\E[ZF^{-1}_Y(U)\id_{\{U\not \in C\}} + Z \E[F^{-1}_Y(U)|U\in C]\id_{\{U  \in C\}}]\\
		&=\E[ZF^{-1}_Y(U)\id_{\{U\not \in C\}}]+\E[\E[ZF^{-1}_Y(U)|U\in C]\id_{\{U  \in C\}}]\\
		&=\E[ZF^{-1}_Y(U)\id_{\{U\not \in C\}}]+\E[ZF^{-1}_Y(U)|U\in C]\P(U\in C)\\
		&=\E[ZF^{-1}_Y(U)\id_{\{U\not \in C\}}]+\E[ZF^{-1}_Y(U)\id_{\{U\in C\}}]=\E[ZF^{-1}_Y(U)].
		\end{aligned}$$
		It follows that $\E[Y|X]=\E[F^{-1}_Y(U)|X]=X$, $\P$-almost surely. If a set of distributions, $\M$, is closed under conditional expectation and $F_Y\in\mathcal{M}$, then $F_{\E[Y|X]}\in\M$, which implies that $F^C_Y=F_X\in\M$. Thus $\M$ is also closed under concentration for all intervals.
%
	\end{proof}

		\begin{proof}[Proof of Proposition \ref{lem:lem-4}]
			(i) Suppose that $\mathcal M$ is closed under concentration for all intervals and $\I$ is a finite. Using \eqref{eq:quantile}, we can see that $F^{\I}$ is the resulting distribution obtained by  sequentially  applying finitely many $C$-concentrations to $F $ over all $C\in \mathcal I$. We thus have $F^\mathcal{I}\in\M$ for all $F\in\M$.
			
			(ii) Suppose that $\M$ is closed under conditional expectation and $F\in\M$. We define
			$$X=F^{-1}(U)\id_{\{U\not \in  \bigcup_{C\in \mathcal I} C\}}  +  \sum_{ C\in \mathcal I}\E[F^{-1}(U)|U\in C] \id_{\{U\in C \}},$$
			whose left-quantile function is given by \eqref{eq:invG3} according to \eqref{eq:quantile}. Following similar argument to the proof of Proposition \ref{lem:lem-2},  for all $\sigma(X)$-measurable random variables $Z$, we have
			$$\begin{aligned}
			\E[XZ]&=\E[ZF^{-1}(U)\id_{\{U\not \in  \bigcup_{C\in \mathcal I} C\}}  +  \sum_{ C\in \mathcal I}Z\E[F^{-1}(U)|U\in C] \id_{\{U\in C \}}]\\
			&=\E[ZF^{-1}(U)\id_{\{U\not \in  \bigcup_{C\in \mathcal I} C\}}]+\sum_{ C\in \mathcal I}\E[\E[ZF^{-1}(U)|U\in C] \id_{\{U\in C \}}]\\
			&=\E[ZF^{-1}(U)\id_{\{U\not \in  \bigcup_{C\in \mathcal I} C\}}]+\sum_{ C\in \mathcal I}\E[ZF^{-1}(U)\id_{\{U\in C \}}]=\E[ZF^{-1}(U)].
			\end{aligned}$$
		 Thus $\E[F^{-1}(U)|X]=X$, $\P$-almost surely, which implies that $F^{\mathcal{I}}=F_X\in\M$. \qedhere
	\end{proof}

	\subsection{Proofs of results in Section \ref{sec:multi}}
	
	\begin{proof}[Proof of Theorem \ref{th:multi}]
To prove the first statement, according to the proof of Theorem \ref{th:1}, it suffices to show that for all increasing $h\in\H$, $\mathbf{X}\in(\mathcal{L}^1)^n$ and $\mathscr{G}\subset\mathscr{F}$, $\rho_h(\E[f(\mathbf{a},\mathbf{X})|\mathscr{G}])\le\rho_h(f(\mathbf{a},\E[\mathbf{X}|\mathscr{G}]))$, which holds directly by Jensen's inequality and monotonicity of $\rho_h$. The second statement holds by Theorem \ref{th:1}. 
The last statement follows from $\rho_h(\E[f(\mathbf{a},\mathbf{X})|\mathscr{G}]) = \rho_h(f(\mathbf{a},\E[\mathbf{X}|\mathscr{G}]))$ and using Theorem \ref{th:1}.
\end{proof}

	\begin{proof}[Proof of Theorem \ref{prop:RRA}]
(i) For all $\mathbf{X}=(X_1,\dots,X_n)\in(\L^1)^n$, take a comonotonic $\widetilde{\mathbf{X}}=(\widetilde{X}_1,\dots,\widetilde{X}_n)\in(\L^1)^n$ such that $\widetilde{X}_i\laweq X_i$ for all $i=1,\dots,n$. It follows that $\E[g(\mathbf{X})]\le\E[g(\widetilde{\mathbf{X}})]$ for all supermodular functions $g:\R^n\to\R$ due to Theorem 5 of \cite{T80}. By Proposition 2.2.5 of \cite{SCB05}, we have $f(\mathbf{a},\mathbf{X})\le_{\rm icx} f(\mathbf{a},\widetilde{\mathbf{X}})$. Moreover, there exists a standard uniform random variable $U$ such that $\widetilde{X}_i=F^{-1}_{\widetilde{X}_i}(U)$ for all $i=1,\dots,n$ and $f(\mathbf{a},\widetilde{\mathbf{X}})=F^{-1}_{f(\mathbf{a},\widetilde{\mathbf{X}})}(U)$ almost surely \citep{D94}.
Take
$$f(\mathbf{a},\widetilde{\mathbf{X}})^{\I_h}=F^{-1}_{f(\mathbf{a},\widetilde{\mathbf{X}})}(U)\id_{\{U\notin \bigcup_{C\in\I_h}C\}}+\sum_{C\in\I_h}\E[F^{-1}_{f(\mathbf{a},\widetilde{\mathbf{X}})}(U)|U\in C]\id_{\{U\in C\}}\sim F_{f(\mathbf{a},\widetilde{\mathbf{X}})}^{\I_h}.$$
It follows that $f(\mathbf{a},\widetilde{\mathbf{X}})^{\I_h}=\E[f(\mathbf{a},\widetilde{\mathbf{X}})|\mathscr{G}]$, where $\mathscr{G}=\sigma(U\id_{\{U\notin \bigcup_{C\in\I_h}C\}})$. Similarly, $\widetilde{X}_i^{\I_h}=\E[\widetilde{X}_i|\mathscr{G}]$ for all $i=1,\dots,n$, where
$$\widetilde{X}_i^{\I_h}=F^{-1}_{\widetilde{X}_i}(U)\id_{\{U\notin \bigcup_{C\in\I_h}C\}}+\sum_{C\in\I_h}\E[F^{-1}_{\widetilde{X}_i}(U)|U\in C]\id_{\{U\in C\}}\sim F_{\widetilde{X}_i}^{\I_h}.$$
Since $f$ is supermodular and positively homogeneous, we have by Theorem 3 of \cite{MM08} that $f(\mathbf{a},\mathbf{X})$ is concave in $\mathbf{X}$. By Jensen's inequality, we have $$f(\mathbf{a},\widetilde{\mathbf{X}})^{\I_h}=\E[f(\mathbf{a},\widetilde{\mathbf{X}})|\mathscr{G}]\le f(\mathbf{a},\E[\widetilde{\mathbf{X}}|\mathscr{G}])= f(\mathbf{a},\widetilde{X}_1^{\I_h},\dots,\widetilde{X}_n^{\I_h}).$$
Thus we have
$$\begin{aligned}
\rho_{h^*}(f(\mathbf{a},\mathbf{X}))\le \rho_{h^*}(f(\mathbf{a},\widetilde{\mathbf{X}}))=\rho_h(f(\mathbf{a},\widetilde{\mathbf{X}})^{\I_h})&\le \rho_h(f(\mathbf{a},\widetilde{X}_1^{\I_h},\dots,\widetilde{X}_n^{\I_h}))\\&\le \sup_{F_{\mathbf{Y}}\in\mathcal{D}(F_1,\dots,F_n)}\sup_{F_1\in\mathcal{F}_1,\dots,F_n\in\mathcal{F}_n}\rho_{h}(f(\mathbf{a},\mathbf{Y})),
\end{aligned}$$
where the first inequality follows from Theorem 4.A.3 of \cite{SS07} and Theorem 5 of \cite{WWW20} and the second equality is by the proof of Theorem \ref{th:1}. Combined with the fact that $$\sup_{F_{\mathbf{X}}\in\mathcal{D}(F_1,\dots,F_n)}\sup_{F_1\in\mathcal{F}_1,\dots,F_n\in\mathcal{F}_n}\rho_{h}(f(\mathbf{a},\mathbf{X}))\le\sup_{F_{\mathbf{X}}\in\mathcal{D}(F_1,\dots,F_n)}\sup_{F_1\in\mathcal{F}_1,\dots,F_n\in\mathcal{F}_n}\rho_{h^*}(f(\mathbf{a},\mathbf{X})),$$
we have \eqref{eq:eq-multi} holds.

(ii) Suppose that the supremum of the right-hand side of \eqref{eq:eq-multi} is attained by some $F_1\in\mathcal{F}_1,\dots,F_n\in\mathcal{F}_n$ and $F_{\mathbf{X}}\in\mathcal{D}(F_1,\dots,F_n)$. For comonotonic $(\widetilde{X}_1,\dots,\widetilde{X}_n)$ such that $\widetilde{X}_i\sim F_i$ for all $i=1,\dots,n$, using the argument in (i), $$\rho_{h^*}(f(\mathbf{a},\mathbf{X}))\le\rho_h(f(\mathbf{a},\widetilde{X}_1^{\I_h},\dots,\widetilde{X}_n^{\I_h})),$$
where $(\widetilde{X}_1^{\I_h},\dots,\widetilde{X}_n^{\I_h})$ is comonotonic and $\widetilde{X}_i^{\I_h}\sim F_i^{\I_h}$ for all $i=1,\dots,n$. Similarly to the proof of Theorem \ref{th:1} (iii), since $\rho_h\le\rho_{h^*}$, we have the supremum of the left-hand side of \eqref{eq:eq-multi} is attained by $F_1^{\mathcal {I}_h},\dots,F_n^{\mathcal {I}_h}$ and $(\widetilde{X}_1^{\I_h},\dots,\widetilde{X}_n^{\I_h})$, which also obtain the supremum of the right-hand side of \eqref{eq:eq-multi} since
$$\begin{aligned}
~~~~~~~~~~~~\rho_{h^*}(f(\mathbf{a},\widetilde{X}_1^{\I_h},\dots,\widetilde{X}_n^{\I_h}))&\le\max_{F_{\mathbf{X}}\in\mathcal{D}(F_1,\dots,F_n)}\max_{F_1\in\mathcal{F}_1,\dots,F_n\in\mathcal{F}_n}\rho_{h^*}(f(\mathbf{a},\mathbf{X}))\\
&=\max_{F_{\mathbf{X}}\in\mathcal{D}(F_1,\dots,F_n)}\max_{F_1\in\mathcal{F}_1,\dots,F_n\in\mathcal{F}_n}\rho_{h}(f(\mathbf{a},\mathbf{X}))\\
&=\rho_{h}(f(\mathbf{a},\widetilde{X}_1^{\I_h},\dots,\widetilde{X}_n^{\I_h}))\le\rho_{h^*}(f(\mathbf{a},\widetilde{X}_1^{\I_h},\dots,\widetilde{X}_n^{\I_h})).~~~~~~~~~~~~~\qedhere
\end{aligned}$$
\end{proof}

	\subsection{Proofs of results in Section \ref{sec:two_con} and related lemmas}
\label{appen:tech}

	In the following, we write $q$ as the H\"older conjugate of $p$. The following lemma closely resembles Theorem 3.4 of \cite{LCLW19} with only an additional statement on the uniqueness of the quantile function of the maximizer.

	\begin{lemma}\label{lem:p}
		For   $h\in \mathcal H^*$, $m\in \R$, $v>0$ and $p>1$, we have
		$$
		\sup_{F_Y\in \mathcal M(p,m,v)} \rho_h(Y)  = m h(1)+  v [h]_q,
		$$
		If  $0<[h]_q<\infty$, the above supremum is attained by a random variable $X$ such that $F_X\in \mathcal M(p,m,v)$ with its quantile function uniquely determined by
		\begin{equation}\label{eq:optimal} \VaR_t(X) = m  + v  \phi_h^q(t),~~t\in (0,1)~~\text{a.e.}\end{equation}
		If  $[h]_q=0$, the above maximum value is attained by any random variable $X$ such that $F_X\in \mathcal M(p,m,v)$. 
	\end{lemma}
	\begin{proof} 
		The only statement that is more than  Theorem 3.4 of \cite{LCLW19}  
		is the uniqueness of the quantile function \eqref{eq:optimal}.
		Without loss of generality, assume $m=0$ and $v=1$.
		Using the H\"older inequality
		\begin{align*} 
			\sup_{F_Y\in\mathcal M(p,0,1)}  \int_0^1  h'(t) \VaR_{1-t}(Y) \d t &=\sup_{F_Y\in\mathcal M(p,0,1)} \int_0^1  (h'(t)-c_{h,q}) \VaR_{1-t}(Y) \d t   \\
			&\le 
			\sup_{F_Y\in\mathcal M(p,0,1)}   {\Vert h'-c_{h,q}\Vert _q \left(\int_0^1 |\VaR_{1-t}(Y)|^p \d t \right)^{1/p}}  =[h]_q. 
		\end{align*} 
		The maximum is attained by $F_X$  only if the above inequality is an equality, which is equivalent to 
		that the function $t\mapsto |\VaR_{1-t}(X)|^p$ is a multiple of $|h'-c_{h,q}|^q$.
		Therefore, 
		$$
		\VaR_t(X) =   \frac{|h'(1-t)-c_{h,q}| ^{q}}{h'(1-t)-c_{h,q}} [h]_q^{1-q} =\phi_h^q(t),~~t\in (0,1)~~\mbox{a.e.}
		$$ 
		Hence, the quantile function of $X$ is uniquely determined by \eqref{eq:optimal}.
	\end{proof}
	
	
	\begin{lemma}\label{lem:p2} 
		For all $h\in \mathcal H$ with $h=\hat{h}$, $m\in \R$, $v>0$ and $p>1$, if $[h^*]_q<\infty$, we have
		$$
		\sup_{F_Y\in \mathcal M(p,m,v)} \rho_h(Y) = \sup_{F_Y\in \mathcal M(p,m,v)} \rho_{h^*}(Y)=   m h(1) +  v [h^*]_q,
		$$
		and the above suprema are simultaneously attained by a random variable $X$ such that $F_X\in \mathcal M(p,m,v)$ with  
		\begin{equation} \label{eq:xx1}\VaR_t(X) = m+ v  \phi_{h^*}^q(t),~~t\in (0,1)~a.e.\end{equation}
	\end{lemma}
	\begin{proof}
		The statement directly follows from Theorem \ref{th:1} and Lemma \ref{lem:p}.
	\end{proof}

	\begin{proof}[Proof of Theorem \ref{th:th1}]
	Together with Theorem \ref{th:1}, Lemmas \ref{lem:p} and \ref{lem:p2} give the statement in Theorem \ref{th:th1} on the supremum.
	Arguments for the infimum are symmetric. For instance, noting that
	$(-h)^*=-h_*$, 
	Theorem \ref{th:1} yields
	\begin{align*}
		\inf_{F_Y\in \mathcal M(p,m,v)} \rho_h(Y) &= - \sup_{F_Y\in \mathcal M(p,m,v)} \rho_{-h}(Y)
		\\&=   - \sup_{F_Y\in \mathcal M(p,m,v)} \rho_{(-h)^*}(Y)
		\\& =   - \sup_{F_Y\in \mathcal M(p,m,v)} \rho_{-h_*}(Y)
		=  \inf_{F_Y\in \mathcal M(p,m,v)} \rho_{h_*}(Y). 
	\end{align*}
	We omit the detailed arguments for the infimum in Theorem \ref{th:th1}. 
	\end{proof}

	\begin{proof}[Proof of Proposition \ref{co:finiteness}]
		Note that  $\rho_h\le \rho_{h^*}$, which is  implied by $h\le h^*$ and  \eqref{eq:disrep}.
		By H\"older's inequality, for any $Y\in \L^p$, using \eqref{eq:disrep3}, we have
		\begin{align*} 
			\int_0^1  {h^*}'(t) \VaR_{1-t}(Y) \d t &=  
			\int_0^1  ({h^*}'(t)-c_{h^*,q}) \VaR_{1-t}(Y) \d t  + c_{h,q}\E[Y] 
			\\&\le 
			[h^*]_q \Vert Y\Vert _p  +c_{h^*,q}\E[Y]<\infty.
		\end{align*}
		The other half of the statement is analogous.
	\end{proof}
	
	\begin{proof}[Proof of Corollary \ref{prop:single-var}]
		We prove the first half (the suprema). The second half is symmetric to the first half. Theorem \ref{th:th1} and Lemma \ref{lem:p2} give
		$$
		\sup_{F_Y\in \mathcal M(p,m,v)} \VaR_\alpha (Y) =  \sup_{F_Y\in \mathcal M(p,m,v)} \ES_\alpha (Y)  =   m +  v [h^*]_q .
		$$
		By Lemma \ref{lem:p}, the corresponding random variable $Z$ which attains 
		$ \ES_\alpha (Z)  =   m +  v  [h^*]_q $
		has left-quantile function
		$$
		F_Z^{-1}(t)= m + v \phi_{h^*}^q (t) 
		= m + v \frac{ \left|\frac{1}{1-\alpha}\id_{(\alpha,1]} (t) -1\right|^{q}} {\frac{1}{1-\alpha}\id_{(\alpha,1]} (t) -1} [h^*]^{1-q}_q,~~t\in[0,1]~~\text{a.e.}
		$$
		Note that $\phi_{h^*}^q (t)$ only takes two values for $t\ge \alpha$ and $t<\alpha$, respectively. 
		Thus  $Z$ is a bi-atomic random variable, and using $\E[Z]=m$,  we have, for some $k_p>0$,
		$$\p\left(Z= m +   \alpha  k_p \right) = 1-\alpha
		\mbox{~~and~~}
		\p\left(Z= m -  (1-\alpha) k_p \right) = \alpha.$$ 
		We note that the number $k_p$ can be determined from $\E[|Z-m|^p]=v^p$, 
		that is, $$k_p= v \left( \alpha^p (1-\alpha)+ (1-\alpha)^p \alpha\right)^{-1/p},$$
		leading to $$
		\sup_{F_Y\in \M(p,m,v)}\VaR_\alpha(Y)=\sup_{F_Y\in \M(p,m,v)}\ES_\alpha(Y)= m+ v\alpha \left( \alpha^p (1-\alpha)+ (1-\alpha)^p \alpha\right)^{-1/p},
		$$ 
		and thus the desired equalities in the   statement on suprema hold. 
	\end{proof}


%
	\end{document}